\documentclass[11pt]{article} 
\usepackage{amsfonts,amsmath,latexsym,amssymb,mathrsfs,amsthm,comment}
\usepackage{graphicx}
\usepackage{xcolor}
\usepackage{booktabs}

\let\OLDthebibliography\thebibliography
\renewcommand\thebibliography[1]{
  \OLDthebibliography{#1}
  \setlength{\parskip}{1pt}
  \setlength{\itemsep}{0pt plus 0.0ex}
}

\def\numberlikeadb{\global\def\theequation{\thesection.\arabic{equation}}}
\numberlikeadb
\newtheorem{theorem}{Theorem}[section]
\newtheorem{lemma}[theorem]{Lemma}

\newtheorem{remark}[theorem]{Remark}

\usepackage{lscape}
\usepackage{caption}
\usepackage{multirow}
\allowdisplaybreaks
\begin{document}

\title{Asymptotic expansions relating to the distribution of the product of correlated normal random variables
}
\author{Robert E. Gaunt\footnote{Department of Mathematics, The University of Manchester, Oxford Road, Manchester M13 9PL, UK, robert.gaunt@manchester.ac.uk; zixin.ye@postgrad.manchester.ac.uk}\:\, and Zixin Y$\mathrm{e}^{*}$}

\date{} 
\maketitle

\vspace{-5mm}

\begin{abstract} Asymptotic expansions are derived for the tail distribution of the product of two correlated normal random variables with non-zero means and arbitrary variances, and more generally the sum of independent copies of such random variables. Asymptotic approximations are also given for the quantile function. Numerical results are given to test the performance of the asymptotic approximations.
\end{abstract}



\section{Introduction}

Let $(X, Y)$ be a bivariate normal random vector with mean vector $(\mu_X,\mu_Y)$, variances $(\sigma_X^2,\sigma_Y^2)$ and correlation coefficient $\rho$. Since the work of \cite{craig,wb32} in the 1930's, the distribution of the product $Z=XY$ has received much attention in the statistics literature (see \cite{gaunt22,np16} for an overview), and has found numerous applications in fields such as statistical mediation analysis \cite{mac}, chemical physics \cite{hey} and condensed matter physics \cite{ach}. The sum $S_n=\sum_{i=1}^nZ_i$, where $Z_1,\ldots,Z_n$ are independent copies of $Z$, also has applications in areas such as quantum cosmology \cite{g96}, astrophysics \cite{man} and electrical engineering \cite{ware}.

In practical applications, a key difficulty in working with the distribution of the product $Z$, or more generally the sum $S_n$, is that the exact distributions of these random variables take rather complicated forms. Indeed, it was not until 2016 that \cite{cui} succeeded in deriving an exact formula for the probability density function (PDF) of $Z$ that holds for general $\mu_X,\mu_Y\in\mathbb{R}$, $\sigma_X,\sigma_Y>0$ and $-1<\rho<1$. Before stating their formula, we introduce some notation that will be used throughout this paper to simplify formulas. We denote $r_X=\mu_X/\sigma_X$, $r_Y=\mu_Y/\sigma_Y$, $s=\sigma_X\sigma_Y$ and
\begin{align*}C_n=C_n(r_X,r_Y,\rho,n)=\exp\bigg\{-\frac{n(r_X^2+r_Y^2-2\rho r_X r_Y)}{2(1-\rho^2)}\bigg\}.
\end{align*}
With this notation, the formula of \cite{cui} for the PDF of the product $Z$ reads
  \begin{align} \label{pdfk}
	    f_{Z}(x) &= \frac{C_1}{\pi}\exp\bigg(\frac{\rho x}{s(1-\rho^2)} \bigg) \sum_{k=0}^\infty\sum_{j=0}^{2k}\frac{x^{2k-j}|x|^{j-k}}{(2k)!s^k(1-\rho^2)^{2k+1/2}
					} \notag \\
               &\quad \times \binom{2k}{j}(r_X-\rho r_Y)^j (r_Y-\rho r_X)^{2k-j}K_{j-k} 
							    \bigg(\frac{|x|}{s(1-\rho^2)}\bigg), \quad x\in\mathbb{R},
   \end{align}
where $K_\nu(x)$  is a modified Bessel function of the second kind (which is defined in Appendix \ref{appa}). Exact formulas for the PDF of $S_n$ in the non-zero mean case were only very recently obtained by \cite{gnp24} (see Section \ref{secpdf} for these formulas). To date, exact formulas for the cumulative distribution function (CDF) are only available for the case $\mu_X=\mu_Y=0$, and even in this restricted setting the formula is given as an infinite series involving the modified Lommel function of the first kind and the modified Bessel function of the second kind (see \cite[Corollary 2.3]{g24}). Only in the case $\rho=\mu_X=\mu_Y=0$ is a closed-form formula available for the CDF (see \cite[Section 2.3]{gaunt22}).

Given that the PDF (\ref{pdfk}) takes a rather complicated form, and that exact formulas for the CDF are unavailable
except for restricted parameter constellations, it is natural to ask for asymptotic approximations for the distribution of the product $Z$, and more generally the sum $S_n$. Such studies have been conducted for a number of other important probability distributions; see, for example, \cite{seg1,seg2,quantilepaper,erlang}. Recently, in \cite{gz23} we undertook such a study, obtaining asymptotic approximations for key distributional properties of the distribution of the product $Z$, including deriving the leading-order asymptotic behaviour of the PDF and tail probabilities in the limit $|x|\rightarrow\infty$, from which asymptotic approximations for the quantile function were deduced. The approximations of \cite{gz23} generalise asymptotic approximations for the product of two correlated mean zero normal variates (see \cite{bc08,gaunt22}); they also complement asymptotic approximations for the PDF \cite{g17,springer} and tail probabilities \cite{l23} for the distribution of the product of $N$ independent mean zero normal random variables. 

In this paper, we extend the results of \cite{gz23} in two natural directions. Firstly, we consider more generally the sum $S_n$, and secondly we are able to derive the entire asymptotic expansion for the PDF of the sum $S_n$ in the limit $|x|\rightarrow\infty$ (Theorem \ref{thm1}). From this asymptotic expansion we deduce the entire asymptotic expansion for the tail probabilities in the limit $|x|\rightarrow\infty$ (Theorem \ref{thm2}), from which we further deduce asymptotic approximations for the quantile function (Theorem \ref{thmq}). Of course setting $n=1$ in our results yields asymptotic expansions for the PDF and tail probabilities of the product $Z$ in the limit $|x|\rightarrow\infty$. Our asymptotic expansions for the PDF and tail probabilities of the product $Z$, and more generally the sum $S_n$, as well as our asymptotic approximations for the quantile function, are simple to implement and computationally efficient. 

The organisation of the paper is as follows.
In Section \ref{secpdf}, we recall recent formulas of \cite{gnp24} for the PDF of the sum $S_n$ that we will apply to derive our asymptotic expansions for the PDF of the sum $S_n$. In Section \ref{sec2.1}, we present our asymptotic expansions for the distribution of the sum $S_n$. Numerical results are given in Section \ref{sec2.2}, which evaluate the performance of these asymptotic approximations. Some preliminary lemmas are stated and proved in Section \ref{sec3}; these results are applied in Section \ref{sec4}, in which we prove the main results of Section \ref{sec2.1}. Finally, we recall some relevant elementary properties of special functions in Appendix \ref{appa}.

\section{Exact formulas for the density}\label{secpdf}

We derive our asymptotic expansions for the PDF of the sum $S_n$ by performing an asymptotic analysis to recent formulas of \cite{gnp24} for the PDF of the sum $S_n$, which we now recall. We remark that even for the purpose of deriving asymptotic expansions for the PDF of the product $Z$ we found it much easier to work with these formulas than with the PDF (\ref{pdfk}) of $Z$; in our recent work \cite{gz23} we performed an asymptotic analysis to the PDF (\ref{pdfk}) and only succeeded in deriving the leading-order behaviour. To state the formulas of \cite{gnp24} in compact form, we require some notation. We denote the sign function by $\mathrm{sgn}(x)$, which is given by $\mathrm{sgn}(x)=1$ for $x>0$, $\mathrm{sgn}(0)=0$, $\mathrm{sgn}(x)=-1$ for $x<0$. We also let $a_{j,k}(x)=k-j$ if $x\geq0$ and $a_{j,k}(x)=j$ if $x<0$. Then, for $x\in\mathbb{R}$,
\begin{align}
f_{S_n}(x)&=\frac{(1-\rho^2)^{n/2-1}C_n}{2^{n-1}s}\exp\bigg(-\frac{|x|}{s(1+\rho\,\mathrm{sgn}(x))}\bigg)\nonumber\\
&\quad\times\sum_{k=0}^\infty\sum_{j=0}^k\frac{(n/8)^k}{k!\Gamma(n/2+a_{j,k}(x))}\binom{k}{j} \bigg(\frac{1+\rho}{1-\rho}\bigg)^j(r_X-r_Y)^{2j} \bigg(\frac{1-\rho}{1+\rho}\bigg)^{k-j}\nonumber\\
&\quad\times (r_X+r_Y)^{2k-2j}\, U\bigg(1-\frac{n}{2}-a_{j,k}(x),2-n-k,\frac{2|x|}{s(1-\rho^2)}\bigg), \label{for1} 
\end{align}
where $U(a,b,x)$ is a confluent hypergeometric function of the second kind (which is defined in Appendix \ref{appa}).
The double infinite series (\ref{for1}) reduces to a single infinite series for the following parameter values.
Suppose that $r_X-r_Y=0$. Then, for $x\in\mathbb{R}$,
\begin{align*}
f_{S_n}(x)&=\frac{(1-\rho^2)^{n/2-1}}{2^{n-1}s}\exp\bigg(-\frac{nr_X^2}{1+\rho}-\frac{|x|}{s(1+\rho\,\mathrm{sgn}(x))}\bigg)\sum_{k=0}^\infty\frac{(n/2)^k}{k!\Gamma(n/2+a_{0,k}(x))}\\
&\quad \times\bigg(\frac{1-\rho}{1+\rho}\bigg)^{k} r_X^{2k}\, U\bigg(1-\frac{n}{2}-a_{0,k}(x),2-n-k,\frac{2|x|}{s(1-\rho^2)}\bigg). 
\end{align*}
Now, suppose that $r_X+r_Y=0$. Then, for $x\in\mathbb{R}$,
\begin{align}
f_{S_n}(x)&=\frac{(1-\rho^2)^{n/2-1}}{2^{n-1}s}\exp\bigg(-\frac{nr_X^2}{1-\rho}-\frac{|x|}{s(1+\rho\,\mathrm{sgn}(x))}\bigg)\sum_{k=0}^\infty\frac{(n/2)^k}{k!\Gamma(n/2+a_{k,k}(x))} \nonumber\\
&\quad\times\bigg(\frac{1+\rho}{1-\rho}\bigg)^{k} r_X^{2k}\, U\bigg(1-\frac{n}{2}-a_{k,k}(x),2-n-k,\frac{2|x|}{s(1-\rho^2)}\bigg). \label{dip}
\end{align}
The double infinite series (\ref{for1}) reduces to a single term in the case $r_X=r_Y=0$: 
\begin{align*}
f_{S_n}(x)=\frac{(1-\rho^2)^{n/2-1}}{2^{n-1}s\Gamma(n/2)}\exp\bigg(-\frac{|x|}{s(1+\rho\,\mathrm{sgn}(x))}\bigg)U\bigg(1-\frac{n}{2},2-n,\frac{2|x|}{s(1-\rho^2)}\bigg), \quad x\in\mathbb{R}.
\end{align*}
On applying the relation (\ref{uk}) followed by the identity (\ref{par}) we see that this formula can be expressed in terms of the modified Bessel function of the second kind: 
\begin{align}\label{0000}
f_{S_n}(x)=\frac{(n/2)^{(n-1)/2}|x|^{(n-1)/2}}{s^{(n+1)/2}\sqrt{\pi(1-\rho^2)}\Gamma(n/2)}\exp\bigg(\frac{\rho  x}{s(1-\rho^2)} \bigg)K_{\frac{n-1}{2}}\bigg(\frac{ |x|}{s(1-\rho^2)}\bigg), \quad x\in\mathbb{R},  
\end{align}
a formula which was derived independently by \cite{gaunt prod, man, np16,wb32}.

Integral representations for the PDF of $S_n$ were also obtained in \cite{gnp24}. Suppose that $|r_X|\not=|r_Y|$. Then, for $x>0$,
 \begin{align} \label{Watson2}
	    f_{S_n}(x) &= D_1(x) \int_0^\infty (t(t+1))^{(n-2)/4} \exp\bigg(-\frac{2x t}{s(1-\rho^2)}\bigg)  \notag \\ 
						 &\quad \times I_{n/2-1}\bigg(\frac{|r_X-r_Y|\sqrt{nxt}}{(1-\rho)\sqrt{s}}\bigg) I_{n/2-1}\bigg(\frac{|r_X+r_Y|\sqrt{nx(1+t)}}{(1+\rho)\sqrt{s}}\bigg)	\,{\rm d} t,
	 \end{align}	
where $I_\nu(x)$ is a modified Bessel function of the first kind (defined in Appendix \ref{appa}) and
   \begin{align*}
	    D_1(x) &= \frac{C_n}{s(1-\rho^2)}\bigg(\frac{n}{4}\bigg)^{1-n/2}|r_X^2-r_Y^2|^{1-n/2}\bigg(\frac{x}{s}\bigg)^{n/2}  \exp\bigg(
			      - \frac{x}{s(1+\rho)}\bigg)\,.
	 \end{align*}
Integral representations in the cases $r_X-r_X=0$ and $r_X+r_Y=0$ are also given in \cite{gnp24}, and can be deduced from the general formula
(\ref{Watson2}) by taking the limits $r_X-r_X\rightarrow0$ and $r_X+r_Y\rightarrow0$, respectively, using the limiting form (\ref{Itend0}). Suppose $r_X-r_Y=0$ and $r_X\not=0$. Then, for $x>0$,
   \begin{align}
	    f_{S_n}(x) &= D_2(x) \int_0^\infty (t^2(t+1))^{(n-2)/4} \exp\bigg(-\frac{2x t}{s(1-\rho^2)}\bigg)I_{n/2-1}\bigg(\frac{2|r_X|\sqrt{nx(1+t)}}{(1+\rho)\sqrt{s}}\bigg)	\,{\rm d} t,\label{222}
	 \end{align}	
where  
   \begin{align*}
	    D_2(x) &= \frac{n^{(2-n)/4}|r_X|^{1-n/2}}{s(1-\rho)^{n/2}(1+\rho)\Gamma(n/2)}\bigg(\frac{x}{s}\bigg)^{(3n-2)/4}\exp\bigg(-\frac{nr_X^2}{1+\rho}
			      - \frac{x}{s(1+\rho)}\bigg)\,.
	 \end{align*}
Now suppose that $r_X+r_Y=0$ and $r_X\not=0$. Then, for $x>0$,
   \begin{align}
	    f_{S_n}(x) &= D_3(x) \int_0^\infty ((t+1)^2t)^{(n-2)/4} \exp\bigg(-\frac{2x t}{s(1-\rho^2)}\bigg) I_{n/2-1}\bigg(\frac{2|r_X|\sqrt{nxt}}{(1-\rho)\sqrt{s}}\bigg)	\,{\rm d} t,\label{333}
	 \end{align}	
where  
   \begin{align*}
	    D_3(x) &= \frac{n^{(2-n)/4}|r_X|^{1-n/2}}{s(1+\rho)^{n/2}(1-\rho)\Gamma(n/2)}\bigg(\frac{x}{s}\bigg)^{(3n-2)/4} \exp\bigg(-\frac{nr_X^2}{1-\rho}
			      - \frac{x}{s(1+\rho)}\bigg)\,.
	 \end{align*} 
Formulas for the PDF $f_{S_n}(x)$ for $x<0$ are obtained by replacing $(x,\rho,r_Y)$ by $(-x,-\rho,-r_Y)$ in equations (\ref{Watson2}), (\ref{222}) and (\ref{333}).
  


\section{Main results: asymptotic expansions for the distribution
}\label{sec2.1}

In Theorem \ref{thm1}, we state asymptotic expansions for the PDF of the sum $S_n$. The formulas are expressed in terms of the confluent hypergeometric function of the first kind $M(a,b,x)$, which is defined in Appendix \ref{appa}. In order to state the theorem, we also need to introduce some notation.
For $a,b\in\mathbb{R}$ and $i,j=0,1,2,\ldots$, we let $g_{i,j}(a,b)$ be the coefficient of $u^iy^{j/2}$ in the Puiseux series expansion of $(1+uy)^a\exp(by^{-1/2}(\sqrt{1+uy}-1))$ about $y=0$, that is $g_{i,j}(a,b)$ is uniquely determined by the expansion
\begin{equation}\label{powerdef}
(1+uy)^a\exp(by^{-1/2}(\sqrt{1+uy}-1))=\sum_{j=0}^\infty\bigg(\sum_{i=\lceil j/2\rceil}^j g_{i,j}(a,b)u^i\bigg)y^{j/2}
\end{equation}
(for further details on Puiseux series see, for example, \cite[p.\ 91]{dav}).  In particular,
\begin{align}\label{gggg}
g_{0,0}(a,b)=1, \quad g_{1,1}(a,b)=\frac{b}{2}, \quad g_{1,2}(a,b)=a, \quad g_{2,2}(a,b)=\frac{b^2}{8}.
\end{align}
We will let $f$ denote the PDF of the sum $S_n$. The ascending factorial (also known as the Pochhammer symbol \cite[Section 5.2 (iii)]{olver}) is given by $(v)_j=\Gamma(v+j)/\Gamma(v)=v(v+1)\cdots(v+j-1)$, and we denote the ceiling function by $\lceil x\rceil$, which is given by $\lceil x\rceil=\min\{m\in\mathbb{Z}\,:\,m\geq x\}$.  With this notation, we can now state our theorem.

\begin{theorem}\label{thm1} Let $\mu_X,\mu_Y\in\mathbb{R}$, $\sigma_X,\sigma_Y>0$, $-1<\rho<1$ and $n\geq1$. 

\vspace{2mm}

\noindent 1. Suppose that $r_X+r_Y\not=0$. Then, as $x\rightarrow\infty$,
\begin{align}f(x)&\sim\frac{s^{-(n+1)/4}C_n}{2\sqrt{2\pi}}\bigg(\frac{1+\rho}{|r_X+r_Y|\sqrt{n}}\bigg)^{(n-1)/2}\exp\bigg(\frac{n}{8}\bigg(\frac{1+\rho}{1-\rho}\bigg)(r_X-r_Y)^2\bigg)\nonumber \\
&\quad\times x^{(n-3)/4}\exp\bigg(\frac{|r_X+r_Y|}{1+\rho}\sqrt{\frac{nx}{s}}-\frac{x}{s(1+\rho)}\bigg)\sum_{\ell=0}^\infty c_\ell(r_X,r_Y,\rho,n)\frac{s^{\ell/2}}{x^{\ell/2}}, \label{expansion1}
\end{align}
where 
$c_0(r_X,r_Y,\rho,n)=1$ 
and, for $\ell\geq1$,
\begin{align*}
c_\ell(r_X,r_Y,\rho,n)
=\exp\bigg(-\frac{n}{8}\bigg(\frac{1+\rho}{1-\rho}\bigg)(r_X-r_Y)^2\bigg)\sum_{j=0}^\ell h_{j,\ell}(r_X,r_Y,\rho,n), 
\end{align*}
with, for $j\geq0$ and $\ell\geq1$,
\begin{align*}
h_{j,\ell}(r_X,r_Y,\rho,n)&=\frac{((3-n)/2)_{\ell-j}((n-1)/2)_{\ell-j}}{(\ell-j)!2^{\ell-j}n^{(\ell-j)/2}}\bigg(\frac{1+\rho}{|r_X+r_Y|}\bigg)^{\ell-j}\sum_{i=\lceil j/2\rceil}^j\bigg(\frac{n}{2}\bigg)_i\bigg(\frac{1-\rho^2}{2}\bigg)^i\\
&\quad\times M\bigg(\frac{n}{2}+i,\frac{n}{2},\frac{n}{8}\bigg(\frac{1+\rho}{1-\rho}\bigg)(r_X-r_Y)^2\bigg) g_{i,j}\bigg(\frac{n-3}{4}-\frac{\ell-j}{2},\frac{|r_X+r_Y|}{1+\rho}\sqrt{n}\bigg).
\end{align*}

\noindent 2. Suppose that $r_X+r_Y=0$. Then, as $x\rightarrow\infty$,
\begin{align}f(x)\sim \frac{x^{n/2-1}}{(2s)^{n/2}\Gamma(n/2)}\exp\bigg(-\frac{nr_X^2}{2}-\frac{x}{s(1+\rho)}\bigg)\sum_{k=0}^\infty d_k(r_X,\rho,n)\frac{s^k}{x^k}, \label{expan}
\end{align}  
where $d_0(r_X,\rho,n)=1$ and, for $k\geq1$,
\begin{align}
d_k(r_X,\rho,n)&=(-1)^k\frac{(1-n/2)_k(n/2)_k}{k!}\bigg(\frac{1-\rho^2}{2}\bigg)^k\exp\bigg(-\frac{n}{2}\bigg(\frac{1+\rho}{1-\rho}\bigg)r_X^2\bigg)\nonumber\\
&\quad\times M\bigg(\frac{n}{2}+k,\frac{n}{2},\frac{n}{2}\bigg(\frac{1+\rho}{1-\rho}\bigg)r_X^2\bigg) \label{rep1} \\
&=(-1)^k\frac{(1-n/2)_k(n/2)_k}{k!}\bigg(\frac{1-\rho^2}{2}\bigg)^k\sum_{j=0}^k\frac{(n/2)^j}{(n/2)_j}\binom{k}{j}\bigg(\frac{1+\rho}{1-\rho}\bigg)^jr_X^{2j}. \label{rep2}
\end{align}

\noindent 3. Suppose that $r_X-r_Y\not=0$. Then, as $x\rightarrow-\infty$,
\begin{align}f(x)&\sim\frac{s^{-(n+1)/4}C_n}{2\sqrt{2\pi}}\bigg(\frac{1-\rho}{|r_X-r_Y|\sqrt{n}}\bigg)^{(n-1)/2}\exp\bigg(\frac{n}{8}\bigg(\frac{1-\rho}{1+\rho}\bigg)(r_X+r_Y)^2\bigg)\nonumber \\
&\quad\times |x|^{(n-3)/4}\exp\bigg(\frac{|r_X-r_Y|}{1-\rho}\sqrt{\frac{n|x|}{s}}+\frac{x}{s(1-\rho)}\bigg)\sum_{\ell=0}^\infty c_\ell(r_X,-r_Y,-\rho,n)\frac{s^{\ell/2}}{|x|^{\ell/2}} \label{for2}.
\end{align} 
\end{theorem}

\noindent 4. Suppose that $r_X-r_Y=0$. Then, as $x\rightarrow-\infty$,
\begin{align}f(x)\sim \frac{|x|^{n/2-1}}{(2s)^{n/2}\Gamma(n/2)}\exp\bigg(-\frac{nr_X^2}{2}+\frac{x}{s(1-\rho)}\bigg)\sum_{k=0}^\infty d_k(r_X,-\rho,n)\frac{s^k}{|x|^k}. \label{part4}
\end{align}

\begin{remark} By making use of the elementary form (\ref{mspecial}) of the confluent hypergeometric function of the first kind and the explicit formulas given in (\ref{gggg}) for the coefficients $g_{i,j}(a,b)$ we obtain the following explicit formulas for the coefficients $c_k=c_k(r_X,r_Y,\rho,n)$, $k=1,2$: 
\begin{align}
c_1&=\frac{n^{3/2}}{8}|r_X+r_Y|(1-\rho)\bigg\{1+\frac{1}{4}\bigg(\frac{1+\rho}{1-\rho}\bigg)(r_X-r_Y)^2\bigg\}-\frac{(n-1)(n-3)}{8\sqrt{n}}\bigg(\frac{1+\rho}{|r_X+r_Y|}\bigg),  \label{c1} \\
c_2&=\frac{1}{128}n^2(n+2)(1-\rho^2)^2\frac{(r_X+r_Y)^2}{(1+\rho)^2}\bigg\{1+\frac{1}{2}\bigg(\frac{1+\rho}{1-\rho}\bigg)(r_X-r_Y)^2 \nonumber \\
&\quad+\frac{n}{16(n+2)}\bigg(\frac{1+\rho}{1-\rho}\bigg)^2(r_X-r_Y)^4\bigg\}+\frac{1}{16}n(n-3)(1-\rho^2)\bigg\{1+\frac{1}{4}\bigg(\frac{1+\rho}{1-\rho}\bigg)(r_X-r_Y)^2\bigg\} \nonumber \\
&\quad-\frac{1}{64}n(n-1)(n-3)(1-\rho^2)\bigg\{1+\frac{1}{4}\bigg(\frac{1+\rho}{1-\rho}\bigg)(r_X-r_Y)^2\bigg\} \nonumber \\
&\quad+\frac{(n+1)(n-1)(n-3)(n-5)}{128n}\bigg(\frac{1+\rho}{r_X+r_Y}\bigg)^2.  \label{c2}
\end{align} 
The following expressions for the coefficients $d_k=d_k(r_X,\rho,n)$, $k=1,2$, are immediate from the representation (\ref{rep2}) of the coefficients:
\begin{align*}
d_1&=\frac{1}{8}n(n-2)(1-\rho^2)\bigg\{1+\bigg(\frac{1+\rho}{1-\rho}\bigg)r_X^2\bigg\}, \\
d_2&=\frac{1}{128}(n+2)n(n-2)(n-4) (1-\rho^2)^2\bigg\{1+2\bigg(\frac{1+\rho}{1-\rho}\bigg)r_X^2+\frac{n}{n+2}\bigg(\frac{1+\rho}{1-\rho}\bigg)^2r_X^4\bigg\}.
\end{align*}  

Observe that for $n=2m$, we have $d_k(r_X,\rho,n)=0$ for all $k\geq m$. Thus, in the appropriate special cases $r_X\pm r_Y=0$, the PDF of $S_{2m}$, $m\geq1$, is given exactly by one of the asymptotic expansions (\ref{rep2}) and (\ref{part4}), truncated at the $(m+1)$-th term, up to a remainder term that is of smaller asymptotic order than $O(|x|^{-\ell})$ for any $\ell>0$.
\end{remark}

\begin{remark} In the case $\mu_X=\mu_Y=0$ (so that $r_X=r_Y=0$), the asymptotic expansions (\ref{expan}) and (\ref{part4}) simplify to
\begin{align*}
f(x)\sim \frac{x^{n/2-1}}{(2s)^{n/2}\Gamma(n/2)}\exp\bigg(\!-\frac{x}{s(1+\rho)}\bigg)\sum_{k=0}^\infty d_k(0,\rho,n)\frac{s^k}{x^k}, \quad x\rightarrow\infty,  \end{align*}
and
\begin{align*}
f(x)\sim \frac{|x|^{n/2-1}}{(2s)^{n/2}\Gamma(n/2)}\exp\bigg(\frac{x}{s(1-\rho)}\bigg)\sum_{k=0}^\infty d_k(0,-\rho,n)\frac{s^k}{|x|^k}, \quad x\rightarrow-\infty,  
\end{align*}
where $d_0(0,\rho,n)=1$ and, for $k\geq1$,
\begin{equation*}
d_k(0,\rho,n)=(-1)^k\frac{(1-n/2)_k(n/2)_k}{k!}\bigg(\frac{1-\rho^2}{2}\bigg)^k.
\end{equation*}
These asymptotic expansions can alternatively be derived very easily by simply applying the asymptotic expansion (\ref{Ktendinfinity}) for the modified Bessel function of the second kind to the formula (\ref{0000}) for the PDF of the sum $S_n$ in the case $\mu_X=\mu_Y=0$.
\end{remark}

\begin{remark} 1. The asymptotic expansion (\ref{for2}) can be obtained from the asymptotic expansion (\ref{expansion1}) by replacing $(r_Y,\rho,x)$ by $(-r_Y,-\rho,-x)$, whilst the asymptotic expansion (\ref{part4}) can be obtained from the asymptotic expansion (\ref{expan}) by replacing $(\rho,x)$ by $(-\rho,-x)$. The asymptotic expansions (\ref{part43}) and (\ref{part44}) of Theorem \ref{thm2} below for the tail probabilities of $S_n$ can be obtained from (\ref{part41}) and (\ref{expan0}), respectively, using the same procedure. This is a consequence of the fact that $S_n$ is a sum of $n$ independent copies of the product $Z$, and that $Z$ is equal in distribution to the product $-X'Y'$, where $(X', Y')$ follows the bivariate normal distribution with means $(\mu_X,-\mu_Y)$, variances $(\sigma_X^2,\sigma_Y^2)$ and correlation coefficient $-\rho$. 

\vspace{2mm}

\noindent 2. From the asymptotic expansions (\ref{expansion1}) and (\ref{expan}), we see that there is a transition in the leading-order asymptotic behaviour of the PDF as $x\rightarrow\infty$ based on whether $r_X+r_Y$ is non-zero or zero. If $r_X+r_Y=0$, we have that $f(x)\sim Ax^{n/2-1}\mathrm{e}^{-ax}$, as $x\rightarrow\infty$, for constants $A,a>0$. On the other hand, if $r_X+r_Y\not=0$, then we have that $f(x)\sim Ax^{(n-3)/4}\mathrm{e}^{-ax+b\sqrt{x}}$, as $x\rightarrow\infty$, for some $A,a,b>0$; thus, the tails get heavier as the quantity $|r_X+r_Y|$ increases. Similar comments apply to the asymptotic expansions (\ref{for2})--(\ref{part4}), as well as the asymptotic expansions of Theorem \ref{thm2} below.
\end{remark}

In the following theorem, we provide asymptotic expansions for the tail probabilities of the sum $S_n$. We denote $F(x)=\mathbb{P}(S_n\leq x)$ and $\bar{F}(x)=1-F(x)=\mathbb{P}(S_n>x)$.

\begin{theorem}\label{thm2}Let $\mu_X,\mu_Y\in\mathbb{R}$, $\sigma_X,\sigma_Y>0$, $-1<\rho<1$ and $n\geq1$. 

\vspace{2mm}

\noindent 1. Suppose that $r_X+r_Y\not=0$. Then, as $x\rightarrow\infty$,
\begin{align}\bar{F}(x)&\sim\frac{(1+\rho)C_n}{2\sqrt{2\pi }}\bigg(\frac{1+\rho}{|r_X+r_Y|\sqrt{n}}\bigg)^{(n-1)/2}\exp\bigg(\frac{n}{8}\bigg(\frac{1+\rho}{1-\rho}\bigg)(r_X-r_Y)^2\bigg)\nonumber \\
&\quad\times\bigg(\frac{x}{s}\bigg)^{(n-3)/4}\exp\bigg(\frac{|r_X+r_Y|}{1+\rho}\sqrt{\frac{nx}{s}}-\frac{x}{s(1+\rho)}\bigg)\sum_{p=0}^\infty \gamma_p(r_X,r_Y,\rho,n)\frac{s^{p/2}}{x^{p/2}}, \label{part41}
\end{align}
where $\gamma_0(r_X,r_Y,\rho,n)=1$ and, for $p\geq1$,
\begin{align*}
\gamma_p(r_X,r_Y,\rho,n)&=\sum_{\substack{i,j,k,\ell\geq0 \\ i+2j+k+\ell=p}}(-1)^{j+i} c_\ell(r_X,r_Y,\rho,n) \binom{(n-1)/2-\ell}{k}\binom{(n-3)/2-\ell-k-2j}{i}\\
&\quad\times ((\ell+k)/2-(n-3)/4)_j(1+\rho)^j\bigg(\frac{|r_X+r_Y|\sqrt{n}}{2}\bigg)^{k+i},     
\end{align*}
with $c_\ell(r_X,r_Y,\rho,n)$, $\ell\geq0$, defined as in part 1 of Theorem \ref{thm1}.

\vspace{2mm}

\noindent 2. Suppose that $r_X+r_Y=0$. Then, as $x\rightarrow\infty$,
\begin{align}\bar{F}(x)&\sim \frac{1+\rho}{2^{n/2}\Gamma(n/2)}\bigg(\frac{x}{s}\bigg)^{n/2-1}\exp\bigg(-\frac{nr_X^2}{2}-\frac{x}{s(1+\rho)}\bigg)\sum_{k=0}^\infty \delta_k(r_X,\rho,n)\frac{s^k}{x^k}, \label{expan0}
\end{align}  
where $\delta_0(r_X,\rho,n)=1$ and, for $k\geq1$,
\begin{align*}
\delta_k(r_X,\rho,n)=\sum_{j=0}^k (-1)^{k-j}d_j(r_X,\rho,n)  (k+1-n/2)_{k-j}(1+\rho)^{k-j},  
\end{align*}
with $d_j(r_X,\rho,n)$, $j\geq0$, defined as in part 2 of Theorem \ref{thm1}.

\vspace{2mm}

\noindent 3. Suppose that $r_X-r_Y\not=0$. Then, as $x\rightarrow-\infty$,
\begin{align}F(x)&\sim\frac{(1-\rho)C_n}{2\sqrt{2\pi }}\bigg(\frac{1-\rho}{|r_X-r_Y|\sqrt{n}}\bigg)^{(n-1)/2}\exp\bigg(\frac{n}{8}\bigg(\frac{1-\rho}{1+\rho}\bigg)(r_X+r_Y)^2\bigg)\nonumber \\
&\quad\times\bigg(\frac{|x|}{s}\bigg)^{(n-3)/4}\exp\bigg(\frac{|r_X-r_Y|}{1-\rho}\sqrt{\frac{n|x|}{s}}+\frac{x}{s(1-\rho)}\bigg)\sum_{p=0}^\infty \gamma_p(r_X,-r_Y,-\rho,n)\frac{s^{p/2}}{x^{p/2}}. \label{part43}
\end{align}
\noindent 4. Suppose that $r_X-r_Y=0$. Then, as $x\rightarrow-\infty$,
\begin{align}F(x)&\sim \frac{1-\rho}{2^{n/2}\Gamma(n/2)}\bigg(\frac{|x|}{s}\bigg)^{n/2-1}\exp\bigg(-\frac{nr_X^2}{2}+\frac{x}{s(1-\rho)}\bigg)\sum_{k=0}^\infty \delta_k(r_X,-\rho,n)\frac{s^k}{|x|^k}. \label{part44}
\end{align} 
\end{theorem}

\begin{remark} We have the following explicit formulas for the coefficients $\gamma_k=\gamma_k(r_X,r_Y,\rho,n)$ and $\delta_k=\delta_k(r_X,\rho,n)$ for $k=1,2$:
\begin{align*}
\gamma_1&=c_1+\frac{1}{2}\sqrt{n}|r_X+r_Y|, \\
\gamma_2&=c_2+\frac{c_1}{2}\sqrt{n}|r_X+r_Y|+\frac{1}{4}(n-3)(1+\rho)+\frac{1}{4}n(r_X+r_Y)^2,
\end{align*}
where $c_1$ and $c_2$ are given as in (\ref{c1}) and (\ref{c2}), and
\begin{align*}
\delta_1&=\frac{1}{2}(n-4)(1+\rho)+\frac{1}{8}n(n-2)(1-\rho^2)\bigg\{1+\bigg(\frac{1+\rho}{1-\rho}\bigg)r_X^2\bigg\}, \\
\delta_2&=\frac{1}{4}(n-6)(n-8)(1+\rho)^2+\frac{1}{16}n(n-2)(n-6)(1+\rho)(1-\rho^2)\bigg\{1+\bigg(\frac{1+\rho}{1-\rho}\bigg)r_X^2\bigg\}\\
&\quad+\frac{1}{128}(n+2)n(n-2)(n-4) (1-\rho^2)^2\bigg\{1+2\bigg(\frac{1+\rho}{1-\rho}\bigg)r_X^2+\frac{n}{n+2}\bigg(\frac{1+\rho}{1-\rho}\bigg)^2r_X^4\bigg\}.
\end{align*}
\end{remark}

Let $0<p<1$. In the following theorem, we provide asymptotic approximations for the quantile function $Q(p)=F^{-1}(p)$ of the sum $S_n$. 

\begin{theorem}\label{thmq}Let $\mu_X,\mu_Y\in\mathbb{R}$, $\sigma_X,\sigma_Y>0$, $-1<\rho<1$ and $n\geq1$. 
Let $q=1-p$. 

\vspace{2mm}

\noindent 1. Suppose that $r_X+r_Y\not=0$. Then, as $p\rightarrow1$,
\begin{align}Q(p)&=s(1+\rho)\bigg\{\ln(1/q)+\frac{|r_X+r_Y|\sqrt{n}}{\sqrt{1+\rho}}\sqrt{\ln(1/q)}+\frac{n-3}{4}\ln(\ln(1/q))\nonumber\\
&\quad+\ln\bigg(\frac{(1+\rho)^{(n+1)/4}}{2\sqrt{2\pi}}\bigg(\frac{1+\rho}{|r_X+r_Y|\sqrt{n}}\bigg)^{(n-1)/2}\bigg)+\frac{n}{8}\bigg(\frac{1+\rho}{1-\rho}\bigg)(r_X-r_Y)^2\nonumber\\
&\quad+\frac{n}{4}\frac{(r_X+r_Y)^2}{1+\rho}-\frac{n}{2(1-\rho^2)}\big(r_X^2+r_Y^2-2\rho r_X r_Y\big)\nonumber\\
&\quad+\frac{(n-3)\sqrt{n}}{8}\frac{|r_X+r_Y|}{\sqrt{1+\rho}}\frac{\ln(\ln(1/q))}{\sqrt{\ln(1/q)}}\bigg\}+O\bigg(\frac{1}{\sqrt{\ln(1/q)}}\bigg). \label{q1}
\end{align}
\noindent 2. Suppose that $r_X+r_Y=0$. Then, as $p\rightarrow1$,
\begin{align}Q(p)&=s(1+\rho)\bigg\{\ln(1/q)+\frac{1}{2}(n-2)\ln(\ln(1/q))+\ln\bigg(\frac{(1+\rho)^{n/2}}{2^{n/2}\Gamma(n/2)}\bigg)-\frac{nr_X^2}{2}\bigg\}\nonumber\\
&\quad+O\bigg(\frac{1}{\ln(1/q)}\bigg). \label{q3}
\end{align}
\noindent 3. Suppose that $r_X-r_Y\not=0$. Then, as $p\rightarrow0$,
\begin{align}Q(p)&=-s(1-\rho)\bigg\{\ln(1/p)+\frac{|r_X-r_Y|\sqrt{n}}{\sqrt{1-\rho}}\sqrt{\ln(1/p)}+\frac{n-3}{4}\ln(\ln(1/p))\nonumber\\
&\quad+\ln\bigg(\frac{(1-\rho)^{(n+1)/4}}{2\sqrt{2\pi}}\bigg(\frac{1-\rho}{|r_X-r_Y|\sqrt{n}}\bigg)^{(n-1)/2}\bigg)+\frac{n}{8}\bigg(\frac{1-\rho}{1-\rho}\bigg)(r_X+r_Y)^2\nonumber\\
&\quad+\frac{n}{4}\frac{(r_X-r_Y)^2}{1-\rho}-\frac{n}{2(1-\rho^2)}\big(r_X^2+r_Y^2-2\rho r_X r_Y\big)\nonumber\\
&\quad+\frac{(n-3)\sqrt{n}}{8}\frac{|r_X-r_Y|}{\sqrt{1-\rho}}\frac{\ln(\ln(1/p))}{\sqrt{\ln(1/p)}}\bigg\}+O\bigg(\frac{1}{\sqrt{\ln(1/p)}}\bigg). \label{q2}
\end{align}
\noindent 4. Suppose that $r_X-r_Y=0$. Then, as $p\rightarrow0$,
\begin{align}Q(p)&=-s(1-\rho)\bigg\{\ln(1/p)+\frac{1}{2}(n-2)\ln(\ln(1/p))+\ln\bigg(\frac{(1-\rho)^{n/2}}{2^{n/2}\Gamma(n/2)}\bigg)-\frac{nr_X^2}{2}\bigg\}\nonumber\\
&\quad+O\bigg(\frac{1}{\ln(1/p)}\bigg). \label{q4}
\end{align}
\end{theorem}

\begin{remark}
The asymptotic approximations (\ref{q1}) and (\ref{q2}) provide the first five terms in the asymptotic expansion of $Q(p)$ as $p\rightarrow1$ and $p\rightarrow0$, respectively, whilst the asymptotic approximations (\ref{q3}) and (\ref{q4}) provide the first three terms in the asymptotic expansions. These asymptotic approximations were derived using only the leading-order term in the asymptotic expansions of Theorem \ref{thm2}; the higher order terms in the asymptotic expansions of the quantile function would involve the correction terms in the asymptotic expansions of Theorem \ref{thm2}. We remark that the asymptotic approximations of \cite{gz23} for the quantile function of the product $Z$ included the first four terms in the asymptotic expansions in the cases $r_X+r_Y\not=0$ and $r_X-r_Y\not=0$. 
Numerical experiments (not reported) for the parameter constellations and $p$-values considered in Section \ref{sec2.2} indicate that including the fifth term in the asymptotic expansion typically leads to more accurate approximations. 
\end{remark}

\section{Numerical results}\label{sec2.2}

\begin{table}[h]
  \centering
  \caption{\footnotesize{Relative error in approximating the PDF of the product $Z$ by the leading-order term, with first-order correction, and with second-order correction in (\ref{expansion1}) (when $\mu_X+\mu_Y\not=0$) and (\ref{expan}) (when $\mu_X+\mu_Y=0$). 
  }}
\label{table1}
\footnotesize{
\begin{tabular}{l*{6}{c}}
\hline
& \multicolumn{6}{c}{$x$} \\
\cmidrule(lr){2-7}
$(\mu_X,\mu_Y,\rho)$ & 2.5 & 5 & 7.5 & 10 & 12.5 & 15 \\
\hline
(1,-1,-0.5)	&	4.4E-02	&	2.3E-02	&	1.6E-02	&	1.2E-02	&	9.7E-03	&	8.1E-03	\\
(1,-1,-0.5)	&	-8.4E-03	&	-2.3E-03	&	-1.1E-03	&	-6.2E-04	&	-4.0E-04	&	-2.8E-04	\\
(1,-1,-0.5)	&	2.9E-03	&	4.2E-04	&	1.3E-04	&	5.8E-05	&	3.0E-05	&	1.8E-05	\\
(1,-1,0)	&	8.2E-02	&	4.5E-02	&	3.1E-02	&	2.4E-02	&	1.9E-02	&	1.6E-02	\\
(1,-1,0)	&	-2.6E-02	&	-7.6E-03	&	-3.6E-03	&	-2.1E-03	&	-1.4E-03	&	-9.6E-04	\\
(1,-1,0)	&	1.4E-02	&	2.2E-03	&	7.1E-04	&	3.1E-04	&	1.7E-04	&	9.8E-05	\\
(1,-1,0.5)	&	1.2E-01	&	6.7E-02	&	4.6E-02	&	3.5E-02	&	2.9E-02	&	2.4E-02	\\
(1,-1,0.5)	&	-4.5E-02	&	-1.3E-02	&	-6.1E-03	&	-3.5E-03	&	-2.3E-03	&	-1.6E-03	\\
(1,-1,0.5)	&	2.6E-02	&	3.9E-03	&	1.3E-03	&	5.5E-04	&	2.9E-04	&	1.7E-04	\\ \hline

(1,0,-0.5)	&	-9.0E-02	&	-7.2E-02	&	-6.1E-02	&	-5.5E-02	&	-5.0E-02	&	-4.6E-02	\\
(1,0,-0.5)	&	2.7E-02	&	1.3E-02	&	8.1E-03	&	5.9E-03	&	4.5E-03	&	3.7E-03	\\
(1,0,-0.5)	&	1.2E-02	&	5.4E-03	&	3.1E-03	&	2.1E-03	&	1.5E-03	&	1.2E-03	\\
(1,0,0)	&	-7.2E-02	&	-4.8E-02	&	-4.0E-02	&	-3.5E-02	&	-3.3E-02	&	-3.0E-02	\\
(1,0,0)	&	2.0E-02	&	1.8E-02	&	1.5E-02	&	1.2E-02	&	1.0E-02	&	8.7E-03	\\
(1,0,0)	&	-2.5E-02	&	-4.9E-03	&	-6.0E-04	&	5.7E-04	&	8.7E-04	&	8.9E-04	\\
(1,0,0.5)	&	-1.0E-01	&	-6.0E-02	&	-4.1E-02	&	-3.2E-02	&	-2.6E-02	&	-2.3E-02	\\
(1,0,0.5)	&	-4.2E-02	&	-1.4E-02	&	-3.0E-03	&	1.9E-03	&	4.1E-03	&	5.0E-03	\\
(1,0,0.5)	&	-9.5E-02	&	-4.2E-02	&	-2.2E-02	&	-1.2E-02	&	-7.5E-03	&	-4.7E-03	\\ \hline
(1,1,-0.5)	&	-2.1E-01	&	-1.6E-01	&	-1.3E-01	&	-1.1E-01	&	-1.0E-01	&	-9.4E-02	\\
(1,1,-0.5)	&	-2.3E-02	&	-1.5E-02	&	-1.1E-02	&	-9.0E-03	&	-7.4E-03	&	-6.4E-03	\\
(1,1,-0.5)	&	1.4E-02	&	4.6E-03	&	2.3E-03	&	1.4E-03	&	9.8E-04	&	7.2E-04	\\
(1,1,0)	&	-1.1E-01	&	-8.8E-02	&	-7.6E-02	&	-6.8E-02	&	-6.2E-02	&	-5.7E-02	\\
(1,1,0)	&	3.0E-02	&	1.4E-02	&	8.4E-03	&	5.9E-03	&	4.5E-03	&	3.7E-03	\\
(1,1,0)	&	1.9E-02	&	8.0E-03	&	4.6E-03	&	3.0E-03	&	2.2E-03	&	1.7E-03	\\
(1,1,0.5)	&	-5.4E-02	&	-4.0E-02	&	-3.4E-02	&	-3.1E-02	&	-2.8E-02	&	-2.6E-02	\\
(1,1,0.5)	&	2.1E-02	&	1.4E-02	&	1.0E-02	&	7.8E-03	&	6.3E-03	&	5.2E-03	\\
(1,1,0.5)	&	-6.0E-03	&	6.5E-04	&	1.1E-03	&	1.0E-03	&	8.1E-04	&	6.6E-04	\\ \hline
    \end{tabular}}
\end{table}

\begin{table}[h]
  \centering
  \caption{\footnotesize{Relative error in approximating $\mathbb{P}(Z\geq x)$ by the leading-order term, with first-order correction, and with second-order correction in (\ref{part41}) (when $\mu_X+\mu_Y\not=0$) and (\ref{expan0}) (when $\mu_X+\mu_Y=0$).
  }}
\label{table3}
\footnotesize{
\begin{tabular}{l*{6}{c}}
\hline
& \multicolumn{6}{c}{$x$} \\
\cmidrule(lr){2-7}
$(\mu_X,\mu_Y,\rho)$ & $q_{0.95}$ & $q_{0.975}$ & $q_{0.99}$ & $q_{0.995}$ & $q_{0.999}$ & $q_{0.9999}$ \\
\hline
        (1,-1,-0.5)	&	5.2E-01	&	3.6E-01	&	2.6E-01	&	2.1E-01	&	1.5E-01	&	9.5E-02	\\
(1,-1,-0.5)	&	-2.5E+00	&	-1.3E+00	&	-7.4E-01	&	-5.4E-01	&	-3.2E-01	&	-2.0E-01	\\
(1,-1,-0.5)	&	1.6E+01	&	5.0E+00	&	1.7E+00	&	9.2E-01	&	2.9E-01	&	5.2E-02	\\
(1,-1,0)	&	4.3E-01	&	3.1E-01	&	2.3E-01	&	1.9E-01	&	1.4E-01	&	9.1E-02	\\
(1,-1,0)	&	-1.7E+00	&	-1.0E+00	&	-6.3E-01	&	-4.7E-01	&	-3.0E-01	&	-1.9E-01	\\
(1,-1,0)	&	8.6E+00	&	3.2E+00	&	1.2E+00	&	6.9E-01	&	2.2E-01	&	3.9E-02	\\
(1,-1,0.5)	&	4.0E-01	&	3.0E-01	&	2.2E-01	&	1.8E-01	&	1.3E-01	&	9.2E-02	\\
(1,-1,0.5)	&	-1.5E+00	&	-9.0E-01	&	-5.7E-01	&	-4.4E-01	&	-2.8E-01	&	-1.8E-01	\\
(1,-1,0.5)	&	6.3E+00	&	2.5E+00	&	1.0E+00	&	5.9E-01	&	1.9E-01	&	3.3E-02	\\\hline
(1,0,-0.5)	&	-3.1E-01	&	-2.9E-01	&	-2.8E-01	&	-2.6E-01	&	-2.4E-01	&	-2.2E-01	\\
(1,0,-0.5)	&	1.2E-01	&	8.5E-02	&	6.0E-02	&	4.7E-02	&	3.4E-02	&	2.5E-02	\\
(1,0,-0.5)	&	1.6E-01	&	1.1E-01	&	7.9E-02	&	6.4E-02	&	4.7E-02	&	3.5E-02	\\
(1,0,0)	&	-1.8E-01	&	-1.7E-01	&	-1.6E-01	&	-1.6E-01	&	-1.5E-01	&	-1.2E-01	\\
(1,0,0)	&	1.8E-01	&	1.4E-01	&	1.0E-01	&	8.9E-02	&	7.0E-02	&	6.4E-02	\\
(1,0,0)	&	7.7E-02	&	6.0E-02	&	4.6E-02	&	3.9E-02	&	3.3E-02	&	3.6E-02	\\
(1,0,0.5)	&	-1.3E-01	&	-1.3E-01	&	-1.2E-01	&	-1.2E-01	&	-1.1E-01	&	-9.0E-02	\\
(1,0,0.5)	&	1.6E-01	&	1.3E-01	&	9.9E-02	&	8.6E-02	&	6.7E-02	&	6.2E-02	\\
(1,0,0.5)	&	7.9E-03	&	1.1E-02	&	1.2E-02	&	1.2E-02	&	1.2E-02	&	2.1E-02	\\ \hline
(1,1,-0.5)	&	-5.5E-01	&	-5.2E-01	&	-4.9E-01	&	-4.7E-01	&	-4.4E-01	&	-4.1E-01	\\
(1,1,-0.5)	&	-1.9E-01	&	-1.7E-01	&	-1.5E-01	&	-1.4E-01	&	-1.2E-01	&	-1.1E-01	\\
(1,1,-0.5)	&	1.5E-03	&	-2.9E-03	&	-5.1E-03	&	-5.1E-03	&	-8.8E-03	&	-2.0E-02	\\
(1,1,0)	&	-3.9E-01	&	-3.7E-01	&	-3.4E-01	&	-3.3E-01	&	-3.1E-01	&	-2.8E-01	\\
(1,1,0)	&	-2.2E-02	&	-2.6E-02	&	-2.6E-02	&	-2.6E-02	&	-2.9E-02	&	-2.6E-02	\\
(1,1,0)	&	8.0E-02	&	5.9E-02	&	4.5E-02	&	3.7E-02	&	2.2E-02	&	1.4E-02	\\
(1,1,0.5)	&	-2.9E-01	&	-2.8E-01	&	-2.6E-01	&	-2.5E-01	&	-2.3E-01	&	-2.2E-01	\\
(1,1,0.5)	&	4.6E-02	&	3.1E-02	&	2.2E-02	&	1.2E-02	&	4.0E-03	&	-1.6E-02	\\
(1,1,0.5)	&	8.4E-02	&	6.2E-02	&	4.8E-02	&	3.4E-02	&	2.2E-02	&	-2.6E-03	\\\hline
    \end{tabular}}
\end{table}

\begin{table}[h]
  \centering
  \caption{\footnotesize{Relative error in approximating $\mathbb{P}(S_n\geq x)$ by the asymptotic expansions (\ref{part41}) (when $\mu_X+\mu_Y\not=0$) and (\ref{expan0}) (when $\mu_X+\mu_Y=0$) with second-order correction.
  }}
\label{table4}
\footnotesize{
\begin{tabular}{l*{6}{c}}
\hline
& \multicolumn{6}{c}{$x$} \\
\cmidrule(lr){2-7}
$(\mu_X,\mu_Y,\rho,n)$ & $q_{0.95}$ & $q_{0.975}$ & $q_{0.99}$ & $q_{0.995}$ & $q_{0.999}$ & $q_{0.9999}$ \\
\hline
(1,-1,-0.5,3)	&	N/A	&	1.4E+00	&	1.6E-01	&	1.7E-02	&	-6.3E-02	&	-7.0E-02	\\
(1,-1,-0.5,5)	&	N/A	&	N/A	&	-2.7E-01	&	-2.5E-01	&	-1.8E-01	&	-1.2E-01	\\
(1,-1,-0.5,7)	&	N/A	&	N/A	&	N/A	&	N/A	&	-3.4E-01	&	-1.8E-01	\\
(1,-1,0,3)	&	5.7E-01	&	1.3E-01	&	-1.3E-02	&	-4.8E-02	&	-6.9E-02	&	-7.1E-02	\\
(1,-1,0,5)	&	-2.9E-01	&	-2.5E-01	&	-2.0E-01	&	-1.7E-01	&	-1.3E-01	&	-9.8E-02	\\
(1,-1,0,7)	&	N/A	&	-4.9E-01	&	-3.2E-01	&	-2.5E-01	&	-1.7E-01	&	-1.2E-01	\\
(1,-1,0.5,3)	&	1.6E-01	&	1.9E-02	&	-4.3E-02	&	-6.0E-02	&	-7.0E-02	&	-7.0E-02	\\
(1,-1,0.5,5)	&	-2.4E-01	&	-2.0E-01	&	-1.6E-01	&	-1.4E-01	&	-1.1E-01	&	-8.3E-02	\\
(1,-1,0.5,7)	&	-3.5E-01	&	-2.7E-01	&	-2.1E-01	&	-1.8E-01	&	-1.3E-01	&	-1.0E-01	\\ \hline
(1,0,-0.5,3)	&	-1.4E-01	&	-1.2E-01	&	-1.0E-01	&	-9.0E-02	&	-7.3E-02	&	-5.3E-02	\\
(1,0,-0.5,5)	&	-5.9E-01	&	-5.1E-01	&	-4.3E-01	&	-3.8E-01	&	-3.1E-01	&	-2.4E-01	\\
(1,0,-0.5,7)	&	-8.7E-01	&	-8.0E-01	&	-7.1E-01	&	-6.6E-01	&	-5.6E-01	&	-4.6E-01	\\
(1,0,0,3)	&	9.1E-03	&	3.0E-03	&	-7.3E-04	&	-7.1E-04	&	-4.0E-03	&	5.1E-03	\\
(1,0,0,5)	&	-2.2E-01	&	-1.8E-01	&	-1.5E-01	&	-1.3E-01	&	-1.1E-01	&	-8.8E-02	\\
(1,0,0,7)	&	-4.4E-01	&	-3.8E-01	&	-3.2E-01	&	-2.9E-01	&	-2.4E-01	&	-1.9E-01	\\
(1,0,0.5,3)	&	2.9E-02	&	2.0E-02	&	1.3E-02	&	1.3E-02	&	3.7E-03	&	4.7E-03	\\
(1,0,0.5,5)	&	-8.6E-02	&	-7.1E-02	&	-5.9E-02	&	-5.2E-02	&	-4.1E-02	&	-3.2E-02	\\
(1,0,0.5,7)	&	-2.0E-01	&	-1.7E-01	&	-1.4E-01	&	-1.2E-01	&	-9.4E-02	&	-7.3E-02	\\ \hline
(1,1,-0.5,3)	&	-4.5E-01	&	-4.1E-01	&	-3.6E-01	&	-3.3E-01	&	-2.7E-01	&	-2.2E-01	\\
(1,1,-0.5,5)	&	-7.5E-01	&	-7.0E-01	&	-6.5E-01	&	-6.2E-01	&	-5.5E-01	&	-4.8E-01	\\
(1,1,-0.5,7)	&	-9.0E-01	&	-8.7E-01	&	-8.3E-01	&	-8.1E-01	&	-7.5E-01	&	-6.8E-01	\\
(1,1,0,3)	&	-1.7E-01	&	-1.5E-01	&	-1.3E-01	&	-1.2E-01	&	-9.1E-02	&	-6.7E-02	\\
(1,1,0,5)	&	-4.0E-01	&	-3.6E-01	&	-3.1E-01	&	-2.9E-01	&	-2.4E-01	&	-1.9E-01	\\
(1,1,0,7)	&	-5.7E-01	&	-5.3E-01	&	-4.8E-01	&	-4.5E-01	&	-3.9E-01	&	-3.3E-01	\\
(1,1,0.5,3)	&	-6.8E-02	&	-5.8E-02	&	-4.9E-02	&	-4.4E-02	&	-3.8E-02	&	-1.4E-02	\\
(1,1,0.5,5)	&	-1.9E-01	&	-1.7E-01	&	-1.4E-01	&	-1.3E-01	&	-1.0E-01	&	-1.0E-01	\\
(1,1,0.5,7)	&	-2.9E-01	&	-2.6E-01	&	-2.3E-01	&	-2.1E-01	&	-1.7E-01	&	-1.5E-01	\\ \hline
    \end{tabular}}
\end{table}

\begin{table}[h]
  \centering
\caption{\footnotesize{Relative error in approximating the quantile function $Q(p)$ by the asymptotic approximations (\ref{q1}) (when $\mu_X+\mu_Y\not=0)$ and (\ref{q3}) (when $\mu_X+\mu_Y=0)$.    }}
\label{table5}
\footnotesize{
\begin{tabular}{l*{6}{c}}
\hline
& \multicolumn{6}{c}{$p$} \\
\cmidrule(lr){2-7}
$(\mu_X,\mu_Y,\rho,n)$ &        0.95 &    0.975 & 0.99 &  0.995 &  0.999 & 0.9999 \\
\hline
(1,-1,-0.5,3)	&	N/A	&	1.0E-01	&	2.2E-02	&	9.1E-03	&	9.8E-04	&	-8.0E-04	\\
(1,-1,-0.5,5)	&	N/A	&	N/A	&	2.6E-01	&	2.4E-02	&	-2.0E-02	&	-1.9E-02	\\
(1,-1,-0.5,7)	&	N/A	&	N/A	&	N/A	&	N/A	&	-5.7E-02	&	-5.2E-02	\\
(1,-1,0,3)	&	-7.4E-02	&	-4.4E-02	&	-2.7E-02	&	-2.1E-02	&	-1.3E-02	&	-8.2E-03	\\
(1,-1,0,5)	&	-8.2E-01	&	-3.3E-01	&	-1.8E-01	&	-1.3E-01	&	-7.8E-02	&	-4.7E-02	\\
(1,-1,0,7)	&	N/A	&	-1.2E+00	&	-4.7E-01	&	-3.2E-01	&	-1.8E-01	&	-1.1E-01	\\
(1,-1,0.5,3)	&	-9.6E-02	&	-6.1E-02	&	-3.9E-02	&	-3.0E-02	&	-1.8E-02	&	-1.1E-02	\\
(1,-1,0.5,5)	&	-4.6E-01	&	-2.9E-01	&	-1.9E-01	&	-1.4E-01	&	-9.0E-02	&	-5.5E-02	\\
(1,-1,0.5,7)	&	-9.9E-01	&	-6.0E-01	&	-3.9E-01	&	-3.0E-01	&	-1.9E-01	&	-1.2E-01	\\ \hline
(1,0,-0.5,3)	&	-1.4E-01	&	-1.2E-01	&	-1.1E-01	&	-9.7E-02	&	-8.2E-02	&	-6.7E-02	\\
(1,0,-0.5,5)	&	-2.4E-01	&	-2.0E-01	&	-1.6E-01	&	-1.5E-01	&	-1.2E-01	&	-9.8E-02	\\
(1,0,-0.5,7)	&	-3.3E-01	&	-2.5E-01	&	-2.0E-01	&	-1.8E-01	&	-1.5E-01	&	-1.2E-01	\\
(1,0,0,3)	&	-1.4E-01	&	-1.1E-01	&	-9.2E-02	&	-8.1E-02	&	-6.4E-02	&	-4.9E-02	\\
(1,0,0,5)	&	-2.4E-01	&	-1.9E-01	&	-1.5E-01	&	-1.3E-01	&	-1.0E-01	&	-8.0E-02	\\
(1,0,0,7)	&	-3.3E-01	&	-2.6E-01	&	-2.0E-01	&	-1.8E-01	&	-1.4E-01	&	-1.0E-01	\\
(1,0,0.5,3)	&	-1.2E-01	&	-1.0E-01	&	-8.1E-02	&	-7.0E-02	&	-5.5E-02	&	-4.1E-02	\\
(1,0,0.5,5)	&	-2.1E-01	&	-1.7E-01	&	-1.3E-01	&	-1.2E-01	&	-8.9E-02	&	-6.7E-02	\\
(1,0,0.5,7)	&	-2.8E-01	&	-2.2E-01	&	-1.8E-01	&	-1.5E-01	&	-1.2E-01	&	-8.9E-02	\\ \hline
(1,1,-0.5,3)	&	-3.6E-01	&	-3.2E-01	&	-2.8E-01	&	-2.6E-01	&	-2.3E-01	&	-1.9E-01	\\
(1,1,-0.5,5)	&	-4.6E-01	&	-4.1E-01	&	-3.6E-01	&	-3.4E-01	&	-2.9E-01	&	-2.5E-01	\\
(1,1,-0.5,7)	&	-5.2E-01	&	-4.7E-01	&	-4.1E-01	&	-3.8E-01	&	-3.3E-01	&	-2.9E-01	\\  
(1,1,0,3)	&	-2.7E-01	&	-2.4E-01	&	-2.1E-01	&	-1.9E-01	&	-1.6E-01	&	-1.3E-01	\\
(1,1,0,5)	&	-3.8E-01	&	-3.3E-01	&	-2.9E-01	&	-2.6E-01	&	-2.2E-01	&	-1.8E-01	\\
(1,1,0,7)	&	-4.5E-01	&	-3.9E-01	&	-3.4E-01	&	-3.1E-01	&	-2.6E-01	&	-2.2E-01	\\
(1,1,0.5,3)	&	-2.2E-01	&	-1.9E-01	&	-1.7E-01	&	-1.5E-01	&	-1.3E-01	&	-1.0E-01	\\
(1,1,0.5,5)	&	-3.2E-01	&	-2.8E-01	&	-2.4E-01	&	-2.2E-01	&	-1.8E-01	&	-1.5E-01	\\
(1,1,0.5,7)	&	-3.9E-01	&	-3.3E-01	&	-2.9E-01	&	-2.6E-01	&	-2.2E-01	&	-1.8E-01	\\ \hline
    \end{tabular}}
\end{table}

In this section, we test the performance of the asymptotic approximations given in Section \ref{sec2.1}. In Table \ref{table1}, we report the relative error in approximating the PDF of $Z$ by the asymptotic expansions (\ref{expansion1}) (when $\mu_X+\mu_Y\not=0$) and (\ref{expan}) (when $\mu_X+\mu_Y=0$). 
Throughout this section, we set $\sigma_X=\sigma_Y=1$ (in this case, $r_X=\mu_X$ and $r_Y=\mu_Y)$. We used \emph{Mathematica} to compute the PDF (\ref{pdfk}) by truncation of the infinite series at $k=50$, which we found to be fast to implement and gave accurate results. In Tables \ref{table3} and \ref{table4}, we report the relative error in approximating the tail probability $\mathbb{P}(S_n\geq x)$ by the asymptotic expansions (\ref{part41}) (for $\mu_X+\mu_Y\not=0$) and (\ref{expan0}) (for $\mu_X+\mu_Y=0$) when $x$ is taken to be a quantile $q_p$ for $p$ taking values from 0.95 to 0.9999. In Table \ref{table5}, we report the relative error in approximating the quantile function of the sum $S_n$ by the asymptotic approximations (\ref{q1}) (when $\mu_X+\mu_Y\not=0$) and (\ref{q3}) (when $\mu_X+\mu_Y=0$). In Tables \ref{table1}--\ref{table5}, a negative number means that the approximation is less than the true value. 

The results in Tables \ref{table3}--\ref{table5} were obtained through Monte Carlo simulations using Python. The distribution of $S_n$ was simulated by first generating $2n$ independent standard normal random variables, and then obtaining a realisation of $S_n$ via the distributional relation $S_n=_d\sum_{i=1}^n(U_i+\mu_X)(\rho U_i+\sqrt{1-\rho^2}V_i+\mu_Y)$, where $U_1,\ldots,U_n$ and $V_1,\ldots,V_n$ are independent standard normal variates (recall that $\sigma_X=\sigma_Y=1$ in this section). For given parameters values, $10^{8}$ realisations of the distribution of the sum $S_n$ were generated from which the results in Tables \ref{table3}--\ref{table5} were derived. In order to estimate $Q(p)$, we took the $k$-th largest realisation as the empirical quantile (where $k=\lfloor N(1-p)\rfloor +1$).  For some parameter constellations, the quantiles $Q(0.95)$, $Q(0.975)$, $Q(0.99)$ and $Q(0.995)$ took negative values. In these cases, the asymptotic expansions (\ref{part41}) and (\ref{expan0}) are not valid (these approximations are only valid if $x>0$). The asymptotic approximations (\ref{q1}) and (\ref{q3}) for the quantile function $Q(p)$ were derived from the approximations (\ref{part41}) and (\ref{expan0}), respectively, and therefore these asymptotic approximations are in turn not valid for $p=0.95, 0.975, 0.99,0.995$. In Tables \ref{table3}--\ref{table5}, we denote these instances by N/A.
We report results to 2 significant figures (s.f.). Most entries in the tables are indeed accurate to 2 s.f. (we repeated the simulations several times in order to verify this). However, for $p=0.9999$ some results in Tables \ref{table3}--\ref{table5} were not accurate to 2 s.f.\ since the error from the simulations was non-negligible in comparison to the error from the asymptotic approximation, particularly when a second-order correction was used. This does to some extent underscore the difficulty in estimating key distributional properties of $S_n$ via Monte Carlo simulations for large $p$-values, 
and, overall, the results of Tables \ref{table3}--\ref{table5} give a good representation of the performance of the asymptotic approximations for a range of $p$-values and parameter values.

It can be seen from Table \ref{table1} that, for the parameter constellations and values of $x$ considered, including a first-order correction yields improved approximations over using just the leading-order term in the asymptotic expansions, and including a second-order correction leads to further improved approximations for the PDF of the product $Z$. However, for smaller values of $x$, including first- and second-order corrections can lead to worse approximations. This is to be expected, as if $x$ is too small then one cannot expect the asymptotic expansions of Theorem \ref{thm1} to provide decent approximations. It is in fact quite surprising how accurate the approximations are even for $x=2.5$, with the largest reported relative error being just 12\% (for $(\mu_X,\mu_Y,\rho)=(1,-1,0.5)$ with only the leading-order term being used in this case). 

We also tested the accuracy of the asymptotic approximations of Theorem \ref{thm2} for the tail probabilities for the same parameter constellations and values of $x$ considered in Table \ref{table1}. We obtained similar results (not reported), which is unsurprising given that the asymptotic expansions of Theorem \ref{thm2} were derived using the asymptotic expansions of Theorem \ref{thm1}. In Table \ref{table3}, we instead chose $x$ to be the quantiles $Q(p)$ with $p$-values ranging from $0.95$ to $0.9999$. All quantiles were positive, but in some cases were rather small; for example, for $(\mu_X,\mu_Y,\rho)=(1,-1,-0.5)$, we obtained via Monte Carlo simulations that $Q(0.95)=0.44$ (to 2 s.f.). In such cases, one cannot expect the asymptotic approximations to give reasonable approximations. 

From Table \ref{table4}, we observe that the asymptotic approximations of Theorem \ref{thm2} for the tail probabilities of the sum $S_n$ tend to be less accurate as $n$ increases. In Table \ref{table4}, we only reported results for the second-order correction (which typically gave the most accurate approximations), but we also observed that the accuracy of approximations using only the leading-order term or applying a first-order correction also tends to decrease as $n$ increases. From Table \ref{table5}, we see that a similar story applies for the asymptotic approximations of Theorem \ref{thmq} for the quantile function of the sum $S_n$ with the accuracy of the approximation tending to decrease as $n$ increases.

\section{Preliminary lemmas}\label{sec3}

We will require the following lemmas in the proofs of our main results of Section \ref{sec2.1}.

\begin{lemma}\label{lem1} (i) Let $a>0$ and $b,m\in\mathbb{R}$ and suppose that the function $g:(0,\infty)\rightarrow\mathbb{R}$ has the following asymptotic expansion:
\begin{equation*}
g(x)\sim x^m\mathrm{e}^{-ax+b\sqrt{x}}\sum_{\ell=0}^\infty\frac{u_\ell}{x^{\ell/2}}, \quad x\rightarrow\infty,   
\end{equation*}
where $u_0,u_1,\ldots$ are real-valued constants. Then, as $x\rightarrow\infty$,
\begin{equation}\label{int1}\int_x^\infty g(t)\,\mathrm{d}t\sim\frac{x^m}{a}\mathrm{e}^{-ax+b\sqrt{x}}\sum_{p=0}^\infty\frac{U_p}{x^{p/2}},
\end{equation}
where, for $p\geq0$,
\begin{align*}
U_p=\sum_{\substack{i,j,k,\ell\geq0 \\ i+2j+k+\ell=p}}(-1)^{j+i}u_\ell \binom{2m+1-\ell}{k}\binom{2m-\ell-k-2j}{i}\frac{((\ell+k)/2-m)_j}{a^j}\bigg(\frac{b}{2a}\bigg)^{k+i}.   
\end{align*}
\noindent (ii) Suppose now that $a>0$ and $m\in\mathbb{R}$ and suppose that the function $h:(0,\infty)\rightarrow\mathbb{R}$ has the following asymptotic expansion:
\begin{equation*}
h(x)\sim x^m\mathrm{e}^{-ax}\sum_{j=0}^\infty\frac{v_j}{x^{j}}, \quad x\rightarrow\infty, 
\end{equation*}
where $v_0,v_1,\ldots$ are real-valued constants. Then, as $x\rightarrow\infty$,
\begin{equation}\label{int2}\int_x^\infty h(t)\,\mathrm{d}t\sim\frac{1}{a}x^m\mathrm{e}^{-ax}\sum_{k=0}^\infty\frac{V_k}{x^k},
\end{equation}
where, for $k\geq0$,
\begin{equation*}
V_k=\sum_{j=0}^kv_j\frac{(k-m)_{k-j}}{(-a)^{k-j}}.    
\end{equation*}
\end{lemma}

The following lemma is a straightforward generalisation of Lemma 3.5 of \cite{gz23}. 

\begin{lemma}\label{lemq} Let $a,A,z>0$ and $b,m\in\mathbb{R}$. Let $g:(0,\infty)\rightarrow\mathbb{R}$ be a function such that $g(x)=O(x^{-1/2})$ as $x\rightarrow\infty$.  Consider the equation
\begin{equation}\label{zeqn}Ax^m\mathrm{e}^{-ax+b\sqrt{x}}\big(1+g(x)\big)=z,
\end{equation}
and notice that there is a unique solution $x$ provided $z$ is sufficiently small. Then, as $z\rightarrow0$,
\begin{align}\label{lemqeqn}x&=\frac{1}{a}\ln(1/z)+\frac{b}{a^{3/2}}\sqrt{\ln(1/z)}+\frac{m}{a}\ln(\ln(1/z))+\frac{b^2}{4a^2}+\frac{\ln(A/a^m)}{a}\nonumber\\
&\quad+\frac{bm}{2a^{3/2}}\frac{\ln(\ln(1/z))}{\sqrt{\ln(1/z)}}+O\bigg(\frac{1}{\sqrt{\ln(1/z)}}\bigg).
\end{align}
In the case $b=0$ and $g(x)=O(x^{-1})$ as $x\rightarrow\infty$, the error in the asymptotic approximation (\ref{lemqeqn}) is of the smaller order $O(1/\ln(1/z))$.
\end{lemma}

\begin{remark}
When $b=0$ and $g(x)=0$ for all $x>0$, an exact solution to equation (\ref{zeqn}) of Lemma \ref{lemq} can be given in terms of the Lambert $W$ function; asymptotic expansions for this function are available in the literature (see \cite{c96}).    
\end{remark}

\noindent{\emph{Proof of Lemma \ref{lem1}.}}
(i) We begin by observing that by interchanging the order of integration and summation we have that, as $x\rightarrow\infty$,
\begin{align}
\int_x^\infty g(t)\,\mathrm{d}t\sim\sum_{\ell=0}^\infty u_\ell I_{a,b,m-\ell/2}(x),  \label{nar1}
\end{align}
where, for $q\in\mathbb{R}$,
\begin{equation*}
I_{a,b,q}(x):=\int_x^\infty t^q\mathrm{e}^{-at+b\sqrt{t}}\,\mathrm{d}t.
\end{equation*}

We now perform an asymptotic analysis of the integral $I_{a,b,q}(x)$ in the limit $x\rightarrow\infty$.
By making the change of variable $y=a(\sqrt{t}-b/(2a))^2$ we have that
\begin{align*}
I_{a,b,q}(x)&=\int_x^\infty t^q\exp\bigg(\!-a\bigg(\sqrt{t}-\frac{b}{2a}\bigg)^2+\frac{b^2}{4a}\bigg)\,\mathrm{d}t \\
&=\frac{\mathrm{e}^{b^2/(4a)}}{a^{q+1}}\int_{a(\sqrt{x}-b/(2a))^2}^\infty y^q\bigg(1+\frac{b}{2\sqrt{ay}}\bigg)^{2q+1}\mathrm{e}^{-y}\,\mathrm{d}y.
\end{align*}
Applying the generalised binomial expansion and then interchanging the order of integration and summation gives that, as $x\rightarrow\infty$,
\begin{align}
I_{a,b,m}(x)&\sim\frac{\mathrm{e}^{b^2/(4a)}}{a^{q+1}}\sum_{k=0}^\infty\binom{2q+1}{k}\bigg(\frac{b}{2\sqrt{a}}\bigg)^k\int_{a(\sqrt{x}-b/(2a))^2}^\infty y^{q-k/2}\mathrm{e}^{-y}\,\mathrm{d}y\nonumber\\
&=\frac{\mathrm{e}^{b^2/(4a)}}{a^{q+1}}\sum_{k=0}^\infty\binom{2q+1}{k}\bigg(\frac{b}{2\sqrt{a}}\bigg)^k\Gamma\bigg(q-\frac{k}{2}+1,a\bigg(\sqrt{x}-\frac{b}{2a}\bigg)^2\bigg),\label{nmb}
\end{align}
where $\Gamma(r,x)=\int_x^\infty t^{r-1}\mathrm{e}^{-t}\,\mathrm{d}t$ is the upper incomplete gamma function. Applying the asymptotic expansion for the upper incomplete gamma function 
\begin{equation*}
\Gamma(r,x)\sim x^{r-1}\mathrm{e}^{-x}\sum_{j=0}^\infty\frac{(-1)^j(1-r)_j}{x^j}, \quad x\rightarrow\infty,    
\end{equation*}
(see \cite[equation 8.11.2]{olver}) to (\ref{nmb}) gives that, as $x\rightarrow\infty$,
\begin{align*}
I_{a,b,q}(x)&\sim\frac{x^q}{a}\mathrm{e}^{-ax+b\sqrt{x}}\sum_{k=0}^\infty\sum_{j=0}^\infty\binom{2q+1}{k}\bigg(\frac{b}{2a}\bigg)^k\frac{(k/2-q)_j}{(-a)^j}\bigg(1-\frac{b}{2a\sqrt{x}}\bigg)^{2q-k-2j}\frac{1}{x^{k/2+j}}.
\end{align*} 
Applying the generalised binomial expansion again yields that, as $x\rightarrow\infty$,
\begin{align}\label{nar2}
I_{a,b,q}(x)&\sim\frac{x^q}{a}\mathrm{e}^{-ax+b\sqrt{x}}\sum_{k=0}^\infty\sum_{j=0}^\infty\sum_{i=0}^\infty\binom{2q+1}{k}\binom{2q-k-2j}{i}\bigg(\frac{b}{2a}\bigg)^{k+i}\frac{(k/2-q)_j}{a^j}\frac{(-1)^{j+i}}{x^{k/2+j+i/2}}.
\end{align} 
Combining (\ref{nar1}) and (\ref{nar2}) now yields the asymptotic expansion (\ref{int1}).

\vspace{2mm}

\noindent(ii) We note that, as $x\rightarrow\infty$,
\begin{align*}
\int_x^\infty h(t)\,\mathrm{d}t\sim\sum_{i=0}^\infty v_i I_{a,0,m-i}(x),
\end{align*}
where, from part (i) of the proof, we have that
\begin{equation*}
I_{a,0,q}(x)\sim \frac{x^q}{a}\mathrm{e}^{-ax}\sum_{j=0}^\infty\frac{(-q)_j}{(-a)^j}\frac{1}{x^{j}}, \quad x\rightarrow\infty. 
\end{equation*}
Therefore, as $x\rightarrow\infty$,
\begin{align*}
\int_x^\infty h(t)\,\mathrm{d}t\sim \frac{x^q}{a}\mathrm{e}^{-ax}\sum_{i=0}^\infty\sum_{j=0}^\infty \frac{(i-m)_j}{(-a)^j}\frac{v_i}{x^{i+j}}= \frac{x^q}{a}\mathrm{e}^{-ax}\sum_{k=0}^\infty\sum_{\ell=0}^k v_\ell\frac{(k-m)_{k-\ell}}{(-a)^{k-\ell}}\frac{1}{x^k},  
\end{align*}
which is the desired asymptotic expansion (\ref{int2}).
\hfill $\Box$

\vspace{2mm}

\noindent{\emph{Proof of Lemma \ref{lemq}.}} Set $w=A/z$ and $h(x)=1+g(x)$. From  (\ref{zeqn}) we obtain that
\begin{align}\label{ax1}ax-b\sqrt{x}=\ln(w)+m\ln(x)+\ln(h(x)).
\end{align}
On solving (\ref{ax1}) using the quadratic formula (taking the positive solution) we obtain that
\begin{align}
x&=\frac{1}{4a^2}\Big(2b^2+2b\sqrt{b^2+4a(\ln(w)+m\ln(x)+\ln(h(x))}\nonumber\\
\label{ax2}&\quad+4a(\ln(w)+m\ln(x)+\ln(h(x))\Big).    
\end{align}
It is clear from (\ref{ax1}) that $x=O(\ln(w))$ as $w\rightarrow\infty$, from which it follows that $\ln(x)=O(\ln(\ln(w)))$ as $w\rightarrow\infty$, and $\ln(h(x))=O(1/\sqrt{\ln(w)})$ as $w\rightarrow\infty$ (since $\ln(1+u)=O(u)$ as $u\rightarrow0$). On applying $\sqrt{1-u}=1-u/2+O(u^2)$ as $u\rightarrow0$ 
to (\ref{ax2}) we obtain, after a simplification, that, as $w\rightarrow\infty$,
\begin{align*}
x
=\frac{1}{a}\ln(w)+\frac{b}{a^{3/2}}\sqrt{\ln(w)}+\frac{b^2}{2a^2}+\frac{m}{a}\ln(x)+\frac{bm}{2a^{3/2}}\frac{\ln(x)}{\sqrt{\ln(w)}}+O\bigg(\frac{1}{\sqrt{\ln(w)}}\bigg).
\end{align*}
Recursively applying this approximation now yields that, as $w\rightarrow\infty$,
\begin{align}\label{ax4}x&=\frac{1}{a}\ln(w)+\frac{b}{a^{3/2}}\sqrt{\ln(w)}+\frac{b^2}{2a^2}+\frac{m}{a}\ln\bigg(\frac{1}{a}\ln(w)+r(w)\bigg)\nonumber\\
&\quad+\frac{bm}{2a^{3/2}}\frac{1}{\sqrt{\ln(w)}}\ln\bigg(\frac{1}{a}\ln(w)+r(w)\bigg) +O\bigg(\frac{1}{\sqrt{\ln(w)}}\bigg),
\end{align}
where $r(w)=O(\sqrt{\ln(w)})$ as $w\rightarrow\infty$. The desired asymptotic expansion (\ref{lemqeqn}) now follows from setting $w=A/z$ in (\ref{ax4}) and using the following asymptotic approximations that were given in the proof of Lemma 3.5 of \cite{gz23}: as $w\rightarrow\infty$,
\begin{align*}
\ln\bigg(\frac{1}{a}\ln(w)+r(w)\bigg)
&=\ln(\ln(w))-\ln(a)+O\bigg(\frac{1}{\sqrt{\ln(w)}}\bigg),
\end{align*}
whilst, as $z\rightarrow0$, 
\begin{align*}
\sqrt{\ln(A/z)}
=\sqrt{\ln(1/z)}+O\bigg(\frac{1}{\sqrt{\ln(1/z)}}\bigg), \quad
\ln(\ln(A/z))
=\ln(\ln(1/z))+O\bigg(\frac{1}{\ln(1/z)}\bigg).
\end{align*}

The final assertion that the error in the asymptotic approximation (\ref{lemqeqn}) is of the smaller order $O(1/\ln(1/z))$ in the case $b=0$ and $g(x)=O(x^{-1})$ as $x\rightarrow\infty$ follows from a straightforward analysis of the order of errors in the approximations made in deriving the asymptotic approximation (\ref{lemqeqn}) when $b=0$ in which we use that, since $g(x)=O(x^{-1})$ as $x\rightarrow\infty$, we now have that $\ln(h(x))=O(1/\ln(1/z))$ as $z\rightarrow\infty$.  We omit the details. 
\hfill $\Box$

\section{Proofs of main results}\label{sec4}

\vspace{2mm}

\noindent{\emph{Proof of Theorem \ref{thm1}.}}
1. We first suppose that $|r_X|\not=|r_Y|$, in which case we can apply the integral representation (\ref{Watson2}) of the PDF of $S_n$. To ease notation a little, throughout this proof we set $s=1$; the general case $s>0$ follows from a simple rescaling.
Applying the asymptotic expansion (\ref{Itendinfinity}) to the integral representation (\ref{Watson2}) of the PDF gives that, as $x\rightarrow\infty$,
\begin{align}
f(x)&\sim\int_0^\infty N_1(x)\sum_{k=0}^\infty\frac{\alpha_k}{x^{k/2}} t^{(n-2)/4}(1+t)^{(n-3)/4-k/2}\exp\bigg(-\frac{2xt}{1-\rho^2}\bigg)\nonumber\\
&\quad\times I_{n/2-1}\bigg(|r_X-r_Y|\frac{\sqrt{nxt}}{1-\rho}\bigg)\exp\bigg(|r_X+r_Y|\frac{\sqrt{nx(1+t)}}{1+\rho}\bigg)\,\mathrm{d}t,  \label{near}
\end{align}
where
\begin{align*}
N_1(x)=\frac{1}{\sqrt{2\pi}}\sqrt{\frac{1+\rho}{|r_X+r_Y|}}\frac{D_1(x)}{(nx)^{1/4}}	
\end{align*}
and
\begin{align*}
\alpha_k=(-1)^k\bigg(\frac{1+\rho}{|r_X+r_Y|\sqrt{n}}\bigg)^ka_k(n/2-1),
\end{align*}
and the constants $a_k(\nu)$, $k\geq0$, are defined in (\ref{akv}).
An interchange of the order of integration and summation, now gives that, as $x\rightarrow\infty$,
\begin{align}
f(x)&\sim N_1(x)\sum_{k=0}^\infty\frac{\alpha_k}{x^{k/2}}\int_0^\infty t^{(n-2)/4}(1+t)^{(n-3)/4-k/2}\exp\bigg(-\frac{2xt}{1-\rho^2}\bigg)\nonumber\\
&\quad\times I_{n/2-1}\bigg(|r_X-r_Y|\frac{\sqrt{nxt}}{1-\rho}\bigg)\exp\bigg(|r_X+r_Y|\frac{\sqrt{nx(1+t)}}{1+\rho}\bigg)\,\mathrm{d}t.\nonumber
\end{align}
Here
we have used the classical term-by-term integration technique for asymptotic expansions of integrals, 
which is applicable when the
integrand has a uniform asymptotic expansion in the integration variable \cite{lopez} (which can be seen here from an application of the limiting forms (\ref{Itend0}) and (\ref{Itendinfinity})). (We remark, however, that asymptotic expansions cannot in general be differentiated, as noted by \cite[p.\ 21]{hinch}.) Making the change of variables $u=xt$ now yields that, as $x\rightarrow\infty$,
\begin{align}
f(x)&\sim N_2(x)\sum_{k=0}^\infty\frac{\alpha_k}{x^{k/2}}\int_0^\infty u^{(n-2)/4}\bigg(1+\frac{u}{x}\bigg)^{(n-3)/4-k/2}\exp\bigg(-\frac{2u}{1-\rho^2}\bigg)\nonumber\\
&\quad\times I_{n/2-1}\bigg(|r_X-r_Y|\frac{\sqrt{nu}}{1-\rho}\bigg)\exp\bigg(|r_X+r_Y|\frac{\sqrt{nx}}{1+\rho}\sqrt{1+\frac{u}{x}}\bigg)\,\mathrm{d}u, \label{expexpu}
\end{align}
where 
\begin{equation*}
N_2(x)=\frac{N_1(x)}{x^{(n+2)/4}}.
\end{equation*}
We now manipulate the second exponential function in the integrand of the integral (\ref{expexpu}) to obtain that, as $x\rightarrow\infty$,
\begin{align}
	f(x)&\sim N_3(x)\sum_{k=0}^\infty\frac{\alpha_k}{x^{k/2}}\int_0^\infty u^{(n-2)/4}\bigg(1+\frac{u}{x}\bigg)^{(n-3)/4-k/2}\exp\bigg(-\frac{2u}{1-\rho^2}\bigg)\nonumber\\
	&\quad\times I_{n/2-1}\bigg(|r_X-r_Y|\frac{\sqrt{nu}}{1-\rho}\bigg)\exp\bigg(|r_X+r_Y|\frac{\sqrt{nx}}{1+\rho}\bigg(\sqrt{1+\frac{u}{x}}-1\bigg)\bigg)\,\mathrm{d}u, \nonumber
\end{align}
where
\begin{equation*}
N_3(x)=\exp\bigg(|r_X+r_Y|\frac{\sqrt{nx}}{1+\rho}\bigg)N_2(x).
\end{equation*}
On recalling the definition (\ref{powerdef}) of the constant $g_{i,j}(a,b)$, $i,j\geq0$, as the coefficient of $u^iy^{j/2}$ in the Puiseux series expansion of  $(1+uy)^a\exp(by^{-1/2}(\sqrt{1+uy}-1))$ about $y=0$, and then interchanging the order of integration and summation, we obtain that, as $x\rightarrow\infty$,
\begin{align}
	f(x)&\sim N_3(x)\sum_{k=0}^\infty\sum_{j=0}^\infty \sum_{i=\lceil j/2\rceil}^j \alpha_kg_{i,j}\bigg(\frac{n-3}{4}-\frac{k}{2},\frac{|r_X+r_Y|}{1+\rho}\sqrt{n}\bigg)\frac{1}{x^{(j+k)/2}}\nonumber\\
	&\quad\times \int_0^\infty u^{(n-2)/4+i}\exp\bigg(-\frac{2u}{1-\rho^2}\bigg) I_{n/2-1}\bigg(|r_X-r_Y|\frac{\sqrt{nu}}{1-\rho}\bigg)\,\mathrm{d}u. \label{beach}
\end{align}
Evaluating the integral in (\ref{beach}) using the integral formula (\ref{integral0}) now gives that, as $x\rightarrow\infty$,
\begin{align}
f(x)&\sim N_4(x)\sum_{k=0}^\infty\sum_{j=0}^\infty \frac{\tilde{h}_{j,k}}{x^{(j+k)/2}},\nonumber \label{lasteq}
\end{align}
where 
\begin{align*}
N_4(x)= \frac{1}{2^{n-1}}n^{(n-2)/4}(1+\rho)^{n/2}(1-\rho)|r_X-r_Y|^{n/2-1}N_3(x)   
\end{align*}
and
\begin{align*}
\tilde{h}_{j,k}&=\alpha_k\sum_{i=\lceil j/2\rceil}^j\bigg(\frac{n}{2}\bigg)_i\bigg(\frac{1-\rho^2}{2}\bigg)^i M\bigg(\frac{n}{2}+i,\frac{n}{2},\frac{n}{8}\bigg(\frac{1+\rho}{1-\rho}\bigg)(r_X-
r_Y)^2\bigg) \\
&\quad\times g_{i,j}\bigg(\frac{n-3}{4}-\frac{k}{2},\frac{|r_X+r_Y|}{1+\rho}\sqrt{n}\bigg) ,  
\end{align*}
which, after some simple algebraic manipulations, can be seen to be the asymptotic expansion (\ref{expansion1}).

In order to conclude that the asymptotic expansion (\ref{expansion1}) holds for $r_X+r_Y\not=0$, we must verify that the asymptotic expansion is also valid for $r_X-r_Y=0$ with $r_X\not=0$. The case $r_X-r_Y=0$ with $r_X\not=0$ is dealt with similarly to the case $|r_X|\not=|r_Y|$ except that we make use of the integral representation (\ref{222}) instead of the integral representation (\ref{Watson2}); we omit the details.

\vspace{2mm}

\noindent 2. Suppose that $r_X+r_Y=0$. To simplify notation, we let $z=(n/2)((1+\rho)/(1-\rho))r_X^2$ and
\begin{equation*}
N(x)=\frac{(1-\rho^2)^{n/2-1}}{2^{n-1}\Gamma(n/2)}\exp\bigg(-\frac{nr_X^2}{1-\rho}-\frac{x}{1+\rho}\bigg).   
\end{equation*}
Now, applying the formula (\ref{dip}) for the PDF in the case $r_X+r_Y=0$ followed by the asymptotic expansion (\ref{Utendinfinity}) gives that, as $x\rightarrow\infty$,
\begin{align*}
f(x)&=N(x)\sum_{k=0}^\infty\frac{z^k}{k!}U\bigg(1-\frac{n}{2},2-n-k,\frac{2x}{1-\rho^2}\bigg)\\
&\sim N(x)\sum_{k=0}^\infty\frac{z^k}{k!}\bigg(\frac{2x}{1-\rho^2}\bigg)^{n/2-1}\sum_{\ell=0}^\infty\frac{(1-n/2)_\ell(n/2+k)_\ell}{\ell!}(-1)^\ell\bigg(\frac{1-\rho^2}{2x}\bigg)^\ell.
\end{align*}
Interchanging the order of summation now gives that, as $x\rightarrow\infty$,
\begin{align*}
f(x)&\sim N(x)\bigg(\frac{2x}{1-\rho^2}\bigg)^{n/2-1}\sum_{\ell=0}^\infty(-1)^\ell\bigg(\frac{1-\rho^2}{2x}\bigg)^\ell\frac{(1-n/2)_\ell}{\ell!}\sum_{k=0}^\infty\frac{(n/2+k)_\ell}{k!}z^k\\
&=N(x)\bigg(\frac{2x}{1-\rho^2}\bigg)^{n/2-1}\sum_{\ell=0}^\infty(-1)^\ell\bigg(\frac{1-\rho^2}{2x}\bigg)^\ell\frac{(1-n/2)_\ell(n/2)_\ell}{\ell!}\sum_{k=0}^\infty\frac{(n/2+\ell)_k}{(n/2)_k}\frac{z^k}{k!}\\
&=N(x)\bigg(\frac{2x}{1-\rho^2}\bigg)^{n/2-1}\sum_{\ell=0}^\infty(-1)^\ell\bigg(\frac{1-\rho^2}{2x}\bigg)^\ell\frac{(1-n/2)_\ell(n/2)_\ell}{\ell!} M\bigg(\frac{n}{2}+\ell,\frac{n}{2},z\bigg)\\
&=N(x)\bigg(\frac{2x}{1-\rho^2}\bigg)^{n/2-1}\mathrm{e}^{z}\sum_{\ell=0}^\infty \frac{d_\ell}{x^\ell}=\frac{x^{n/2-1}}{2^{n/2}\Gamma(n/2)}\exp\bigg(-\frac{nr_X^2}{2}-\frac{x}{1+\rho}\bigg)\sum_{\ell=0}^\infty \frac{d_\ell}{x^\ell},
\end{align*}
where the coefficients $d_\ell=d_\ell(r_X,\rho,n)$, $\ell\geq0$, are defined as in (\ref{rep1}). Here we evaluated the sum over the index $k$ by using the hypergeometric series representation of the confluent hypergeometric function of the first kind (combine (\ref{mdef}) and (\ref{gauss})). We have thus derived the asymptotic expansion (\ref{expan}) with representation (\ref{rep1}) for the coefficients $d_k$, $k\geq0$. The second representation (\ref{rep2}) for the coefficients $d_k$, $k\geq0$, follows from the reduction formula (\ref{mspecial}). 

\vspace{2mm}

\noindent 3. Suppose that $r_X-r_Y\not=0$. We now consider the case $x\rightarrow-\infty$. We note that $S_n$ is a sum of $n$ independent copies of the product $Z$ and that $Z=XY=_d-X'Y'$, where $(X', Y')$ follows the bivariate normal distribution with means $(\mu_X,-\mu_Y)$, variances $(\sigma_X^2,\sigma_Y^2)$ and correlation coefficient $-\rho$. We thus obtain the asymptotic expansion (\ref{for2}) by replacing $(x,r_Y,\rho)$ with $(-x,-r_Y,-\rho)$ in the expansion (\ref{expansion1}). 

\vspace{2mm}

\noindent 4. This is similar to part 3 of the proof, except that we now obtain the asymptotic expansion (\ref{part4}) by replacing $(x,\rho)$ by $(-x,-\rho)$ in the asymptotic expansion (\ref{expan}).
\hfill $\Box$

\vspace{3mm}

\noindent{\emph{Proof of Theorem \ref{thm2}.}} We obtain the asymptotic expansion (\ref{part41}) by using that $\bar{F}(x)=\int_x^\infty f(t)\,\mathrm{d}t$ and then applying the asymptotic expansion (\ref{expansion1}) followed by an application of part (i) of Lemma \ref{lem1}. We obtain (\ref{expan0}) similarly, but instead we apply the asymptotic expansion (\ref{expan}) and part (ii) of Lemma \ref{lem1}. We derive (\ref{part43}) and (\ref{part44}) similarly, although we now use that $F(x)=\int_{-\infty}^x f(t)\,\mathrm{d}t$ and apply the asymptotic expansions (\ref{for2}) and (\ref{part4}), instead of (\ref{expansion1}) and (\ref{expan}). \hfill $\Box$

\vspace{3mm}

\noindent{\emph{Proof of Theorem \ref{thmq}.}} 1. Suppose that $r_X+r_Y\not=0$. We begin by recalling that $Q(p)$ solves the equation $\bar{F}(Q(p))=1-p$. On applying the leading-order term in the asymptotic expansion (\ref{part41}), it can be seen that the quantile function $Q(p)$ solves an equation of the form (\ref{zeqn}) with
\begin{align*}
z&=1-p, \quad a=\frac{1}{s(1+\rho)},\quad b=\frac{|r_X+r_Y|\sqrt{n}}{(1+\rho)\sqrt{s}},\quad m=\frac{n-3}{4},\\
A&=\frac{(1+\rho)C_n}{2\sqrt{2\pi }s^{(n-3)/4}}\bigg(\frac{1+\rho}{|r_X+r_Y|\sqrt{n}}\bigg)^{(n-1)/2}\exp\bigg(\frac{n}{8}\bigg(\frac{1+\rho}{1-\rho}\bigg)(r_X-r_Y)^2\bigg).
\end{align*}
We now obtain the desired asymptotic approximation (\ref{q1}) by applying the asymptotic approximation (\ref{lemqeqn}) with these values of $z$, $a$, $b$, $m$ and $A$.

 \vspace{2mm}

 \noindent 2. Suppose that $r_X+r_Y=0$. This is similar to part 1, although in deriving the asymptotic approximation (\ref{q3}) we now use the leading-order term in the asymptotic expansion (\ref{expan0}), in which case
 \begin{align*}
z=1-p,\quad a=\frac{1}{s(1+\rho)}, \quad b=0, \quad m=\frac{n}{2}-1, \quad A=\frac{(1+\rho)}{2^{n/2}s^{n/2-1}\Gamma(n/2)}\mathrm{e}^{-nr_X^2/2}.   
 \end{align*}

\vspace{2mm}

\noindent 3 \& 4. We derive the asymptotic approximations (\ref{q2}) and (\ref{q4}) similarly to the derivations in parts 1 and 2. This time we use the fact that the quantile function $Q(p)$ satisfies $F(Q(p))=p$ and apply the leading-order term in the asymptotic expansions (\ref{part43}) and (\ref{part44}) for the cases $r_X-r_Y\not=0$ and $r_X-r_Y=0$, respectively.  \hfill $\Box$

\section*{Acknowledgements}
RG is funded in part by EPSRC grant EP/Y008650/1 and EPSRC grant UKRI068. ZY is supported by a University of Manchester Research Scholar Award. We would like to thank the reviewers for their helpful comments and suggestions.

\appendix

\section{Special functions}\label{appa}
In this appendix, we recall some elementary properties of the generalized hypergeometric function, the confluent hypergeometric functions of the first and second kind and the modified Bessel functions of the first and second kind.
Unless otherwise stated, these properties can be located in the standard reference \cite{olver}. 

The \emph{generalized hypergeometric function} is defined, for $|x|<1$, by 
\begin{equation}
\label{gauss}
{}_pF_q(a_1,\ldots,a_p; b_1,\ldots,b_q;x)=\sum_{j=0}^\infty\frac{(a_1)_j\cdots(a_p)_j}{(b_1)_j\cdots(b_q)_j}\frac{x^j}{j!},
\end{equation}
and by analytic continuation elsewhere. 
   
The \emph{confluent hypergeometric function of the first kind} can be defined by
\begin{equation}\label{mdef}
M(a,b,x)={}_{1}F_1(a;b;x).    
\end{equation}
The \emph{confluent hypergeometric function of the second kind} can be defined by
\begin{equation*}
U(a,b,x)=\frac{\Gamma(b-1)}{\Gamma(a)}x^{1-b}M(a-b+1,2-b,x)+\frac{\Gamma(1-b)}{\Gamma(a-b+1)}M(a,b,x), \quad b\notin\mathbb{Z},   
\end{equation*}
and
\begin{equation*}
U(a,b,x)=\lim_{\beta\rightarrow b}U(a,\beta,x),\quad b\in\mathbb{Z}.   
\end{equation*}
For integer $m=0,1,2,\ldots$, we have the elementary form:
\begin{equation}
\label{mspecial} M(a+m,a,x)=\mathrm{e}^x\sum_{j=0}^m\binom{m}{j}\frac{x^j}{(a)_j},    
\end{equation}
which can be obtained by applying the relation $M(a,b,x)=\mathrm{e}^xM(b-a,b,-x)$ followed by the series representation of the confluent hypergoemetric function of the first kind (combine (\ref{mdef}) and (\ref{gauss})).
 In particular, we have the following special cases:
\begin{align*}
M(a,a,x)=\mathrm{e}^x, \quad M(a+1,a,x)=\bigg(1+\frac{x}{a}\bigg)\mathrm{e}^x, \quad M(a+2,a,x)=\bigg(1+\frac{2x}{a}+\frac{x^2}{a(a+1)}\bigg)\mathrm{e}^x.
\end{align*}
The confluent hypergeometric function of the second kind has the following asymptotic behaviour:
\begin{align}\label{Utendinfinity} U(a,b,x) \sim x^{-a}\sum_{s=0}^\infty\frac{(a)_s(a-b+1)_s}{s!}(-x)^{-s}, \quad x \rightarrow \infty.
\end{align}

The \emph{modified Bessel function of the first kind} is defined, for $\nu\in\mathbb{R}$ and $x\in\mathbb{R}$, by the power series
\begin{equation}\label{idef}
I_\nu(x)=\sum_{k=0}^\infty \frac{(x/2)^{2k+\nu}}{k!\Gamma(k+\nu+1)}.    
\end{equation}
The \emph{modified Bessel function of the second kind} can be defined, for $\nu\in\mathbb{R}$ and $x>0$, by the integral
\[K_\nu(x) = \frac12 \Big(\frac{x}2\Big)^\nu \int_0^\infty \exp\bigg(-t-\frac{x^2}{4 t}\bigg)\, 
			           \frac{{\rm d} t}{t^{\nu+1}}.
\]
We have the following identity: for $\nu\in\mathbb{R}$ and $x>0$,
\begin{equation}
\label{par} K_{-\nu}(x)=K_{\nu}(x),  
\end{equation}
and the relation
\begin{align}
\label{uk} U(a,2a,2x)&=\frac{1}{\sqrt{\pi}}\mathrm{e}^{x}(2x)^{1/2-a}K_{a-\frac{1}{2}}(x).
\end{align}
The modified Bessel functions possess the following asymptotic behaviour:
\begin{eqnarray}\label{Itend0}I_{\nu} (x) &\sim&\frac{x^\nu}{2^{\nu}\Gamma(\nu+1)}, \quad x\rightarrow0,\: \nu>-1. \\
\label{Itendinfinity} I_{\nu} (x) &\sim& \frac{\mathrm{e}^x}{\sqrt{2\pi x}} \sum_{k=0}^\infty(-1)^k\frac{a_k(\nu)}{x^k}, \quad x \rightarrow \infty,\: \nu\in\mathbb{R}, \\
\label{Ktendinfinity} K_{\nu} (x) &\sim& \sqrt{\frac{\pi}{2x}}\mathrm{e}^{-x} \sum_{k=0}^\infty\frac{a_k(\nu)}{x^k}, \quad x \rightarrow \infty,\: \nu\in\mathbb{R},
\end{eqnarray}
where $a_0(\nu)=1$ and, for $k\geq1$,
\begin{align}\label{akv}a_k(\nu)
=(-1)^k\frac{(1/2-\nu)_k(1/2+\nu)_k}{k!2^k}.
\end{align}
Here the expansion (\ref{Itendinfinity}) is valid for $|\mathrm{ph}(x)|\leq\pi/2-\delta$, whilst the expansion (\ref{Ktendinfinity}) holds for $|\mathrm{ph}(x)|\leq 3\pi/2-\delta$, where $\delta$ denotes an arbitrary small positive constant.

The following definite integral is given in \cite[equation 6.643(2)]{g07}: for $b,\alpha>0$ and $\mu,\nu\in\mathbb{R}$ such that $\mu+\nu+1/2>0$,
\begin{equation}
\label{integral0} \int_0^\infty x^{\mu-1/2}\mathrm{e}^{-\alpha x}I_{2\nu}(2b\sqrt{x})\,\mathrm{d}x=\frac{\Gamma(\mu+\nu+1/2)}{\Gamma(2\nu+1)}\frac{b^{2\nu}}{\alpha^{\mu+\nu+1/2}}M\bigg(\mu+\nu+\frac{1}{2},2\nu+1,\frac{b^2}{\alpha}\bigg).    
\end{equation}
Note that equation 6.643(2) of \cite{g07} is expressed in terms of the Whittaker function of the first kind; we used the standard relation between the Whittaker function of the first kind and the confluent hypergeometric function of the first kind (see \cite[13.14.2]{olver}) to obtain formula (\ref{integral0}). 
Since $M(a,a,x)=\mathrm{e}^x$,
we have the following specialisation of formula (\ref{integral0}): for $b,\alpha,\mu>0$,
\begin{equation}
\label{integral} \int_0^\infty x^{\mu-1/2}\mathrm{e}^{-\alpha x}I_{2\mu-1}(2b\sqrt{x})\,\mathrm{d}x=\frac{b^{2\mu-1}}{\alpha^{2\mu}}\exp\bigg(\frac{b^2}{\alpha}\bigg).    
\end{equation}


\footnotesize

\begin{thebibliography}{99}
\addcontentsline{toc}{section}{References}

\bibitem{ach} Acharya, P., Sengupta, S., Chakraborty, B. and Ramola, K. Athermal fluctuations in disordered crystals. \emph{Phys. Rev. Lett.} $\mathbf{124}$ (2020), 168004.

\bibitem{bc08} Bandi, M. M. and Connaughton, C. Craig's $XY$ distribution and the statistics of Lagrangian power in two-dimensional turbulence. \emph{Phys. Rev. E} $\mathbf{77}$ (2008), 036318.


\bibitem{c96} Corless, R. M., Gonnet, G. H., Hare, D. E., Jeffrey, D. J. and Knuth, D. E.  On the Lambert $W$ function. \emph{Adv. Comput. Math.} $\mathbf{5}$ (1996), 329--359.

\bibitem{craig} Craig, C. C. On the Frequency Function of $xy$. \emph{Ann. Math. Stat.} $\mathbf{7}$ (1936), 1--15.

\bibitem{cui} Cui, G., Yu, X. Iommelli, S. and Kong, L. Exact Distribution for the Product of Two Correlated Gaussian Random Variables.  \emph{IEEE Signal Process. Lett.}  $\mathbf{23}$ (2016), 1662--1666.

\bibitem{dav} Davenport, J. H., Siret, Y. and Tournier, E. \emph{Computer Algebra: Systems and Algorithms for Algebraic Computation,} 2nd ed. San Diego: Academic Press, 1993.

\bibitem{g17} Gaunt, R. E. On Stein's method for products of normal random variables and zero bias couplings. \emph{Bernoulli} $\mathbf{23}$ (2017), 3311--3345.

\bibitem{gaunt prod} Gaunt, R. E. A note on the distribution of the product of zero mean correlated normal random variables. \emph{Stat. Neerl.} $\mathbf{73}$ (2019),  176--179. 

\bibitem{gaunt22} Gaunt, R. E. The basic distributional theory for the product of zero mean correlated normal random variables.  \emph{Stat. Neerl.} $\mathbf{76}$ (2022), 450--470.

\bibitem{g24} Gaunt, R. E. On the cumulative distribution function of the variance-gamma distribution. \emph{B. Aust. Math. Soc.} $\mathbf{110}$ (2024), 389--397.


\bibitem{gnp24} Gaunt, R. E., Nadarajah, S. and Pogany, T. K. Infinite Divisibility of the Product of Two Correlated Normal Random Variables and Exact Distribution of the Sample Mean.  arXiv:2405.10178, 2024.

\bibitem{gz23} Gaunt, R. E. and Ye, Z. Asymptotic approximations for the distribution of the product of correlated normal random variables. \emph{J. Math. Anal. Appl.} $\mathbf{543}$ (2025), Art.\ 128987.

\bibitem{seg1} Gil, A., Segura, J. and Temme, N. M. On the computation and inversion of
the cumulative noncentral beta distribution. \emph{Appl. Math. Comput.} $\mathbf{361}$ (2019), 74--86.

\bibitem{seg2} Gil, A., Segura, J. and Temme, N. M. New asymptotic representations of the noncentral $t$-distribution. \emph{Stud. Appl. Math.} $\mathbf{151}$ (2023), 857--882.

\bibitem{g07} Gradshteyn, I. S. and Ryzhik, I. M.  \emph{Table of Integrals, Series and Products,}  $7$th ed.  Academic Press, 2007.

\bibitem{g96} Grishchuk, L. P. Statistics of the microwave background anisotropies caused by the squeezed cosmological perturbations. \emph{Phys. Rev. D} $\mathbf{53}$ (1996), 6784.


\bibitem{hey} Heyes, D. M., Dini, D. and Smith, E. R. Single trajectory transport coefficients and the energy landscape by molecular dynamics simulations. \emph{J. Chem. Phys.} $\mathbf{152}$, (2020), 194504.

\bibitem{hinch} Hinch, E. J. \emph{Perturbation Methods.} Cambridge Texts in Applied Mathematics. Cambridge University Press, Cambridge, 1991.

\bibitem{l23} Leipus, R., \v{S}iaulys, J., Dirma, M. and Zov\'e, R. On the distribution-tail behaviour of the product of normal random variables. \emph{J. Inequal. Appl.} (2023), no.\ 32.

\bibitem{lopez} L\'opez, J. L. Asymptotic expansions of integrals:
The term-by-term integration method. \emph{J. Comput. Appl. Math.} $\mathbf{102}$ (1999), 181--194. 

\bibitem{mac} MacKinnon, D. P. \emph{Introduction to Statistical Mediation Analysis.} New York, NY: Routledge, 2012.

\bibitem{man} Mangilli, A., Plaszczynski, S. and Tristram, M. Large-scale cosmic microwave background temperature and polarization cross-spectra likelihoods. \emph{Mon. Not. R. Astron. Soc.} $\mathbf{453}$ (2015), 3174--3189.

\bibitem{np16} Nadarajah, S. and Pog\'{a}ny, T. K. On the distribution of the product of correlated normal random variables. \emph{C. R. Acad. Sci. Paris, Ser. I} $\mathbf{354}$ (2016),  201--204.

\bibitem{olver} Olver, F. W. J., Lozier, D. W., Boisvert, R. F. and Clark, C. W.  \emph{NIST Handbook of Mathematical Functions.} Cambridge University Press, 2010.


\bibitem{springer} Springer, M. D. and Thompson, W. E. The distribution of products of Beta, Gamma
and Gaussian random variables. \emph{SIAM J. Appl. Math.} $\mathbf{18}$ (1970), 721--737.

\bibitem{quantilepaper} Stupfler G. and Usseglio-Carleve A. Simple sufficient criteria for second-order extended regular variation in the Gumbel domain of attraction: The case of Weibull-tailed distributions. To appear in \emph{Extremes}, 2025+.

\bibitem{erlang} van Leeuwaarden, J. S. H. and Temme, N. M. Asymptotic inversion of
the Erlang B formula. \emph{SIAM J. Appl. Math.} $\mathbf{70}$ (2009), 1--23.
  
\bibitem{ware} Ware, R. and Lad, F. Approximating the distribution for sums of products of normal variables. Working paper, Department of Mathematics and Statistics, University of Canterbury, 2013.


\bibitem{wb32} Wishart, J. and Bartlett, M. S. The distribution of second order moment statistics in a normal system. \emph{Proc. Camb. Philol. Soc.} $\mathbf{28}$ (1932), 455--459.

\end{thebibliography}

\normalsize
\end{document}